\documentclass[oneside,reqno]{amsart}
\usepackage{fix-cm}
\usepackage{lucacommands}

\title{Determinateness in a classical theory of nonclassical truth}
\author{Luca Castaldo}

\begin{document}
\allowdisplaybreaks
\maketitle
\begin{abstract}
Recent theories that take determinateness as primitive are either nonclassical
theories of nonclassical truth, like Field's \INT, or classical theories of classical truth, like Fujimoto and Halbach's \cd and \cdp. This paper develops a third option: a \textit{classical} theory of \textit{nonclassical} truth. The theory, \kfdcs, extends the Kripke-Feferman theory \kfcs with a primitive determinateness predicate. Unlike \cdp, \kfdcs keeps the iteration of
truth unrestricted and gives up its completeness. 
With or without a consistency axiom, \kfd is proof-theoretically equivalent to the theory \idh2 of twice iterated fixed points, with ordinal $\varphi_{\varphi_{\varepsilon_0}0}0$, thus stronger than \cdp.
It is also shown that replacing full induction with internal induction yields ordinal $\varphi_{\varphi_\omega0}0$, and restricting internal induction to total formulae yields $\varphi_{\varepsilon_0}0$, the ordinal of \kf.

\end{abstract}

\section{Introduction}
Axiomatic theories of truth are severely constrained by the semantic paradoxes. By Tarski's theorem, no consistent theory formulated in classical logic can contain a type-free truth predicate that is fully disquotational. As clearly formulated by Feferman \cite{fef84Toward}, a consistent theory of truth must therefore restrict the language, the underlying logic, or the basic truth-theoretic principles.

A natural response in this context is to try to distinguish between paradoxical and non-paradoxical sentences. There is, the thought goes, a class of `non-paradoxical', `grounded', `meaningful', `determinate' sentences, for which the restrictions forced on us by the semantic paradoxes are no longer necessary: non-paradoxical sentences can obey the principles that fail in general. In this sense, they are \textit{admissible}---as instances of our naive principles.\footnote{
    In the literature, the predicate `healthy' is also used to refer to non-paradoxical sentences following \cite{bac15Can}; here we follow \cite{cgw26}.
}

Admissibility has a long tradition of axiomatic characterizations, many of which take it to be definable in terms of truth.
One of the first examples is due to Reinhardt \cite{rei85Remarks,rei86Remarks}. In his well-known instrumental account of \kfcs---Kripke-Feferman truth with consistency---Reinhardt identifies the property of `significance' induced by \kfcs: the significant sentences, for which restrictions can be lifted, are those that are true or false.
A similar characterization is suggested by Feferman's theory \msf{DT} \cite{fef08Axioms}---Determinate Truth---in which likewise `determinate' sentences are defined as true or false (see also \cite{fef84Toward,fef91Reflecting}). Unlike Reinhardt, however, Feferman takes determinateness to be \textit{conceptually prior} to truth: ``the conditions on \D should be prior to those on \T, that is, determinate meaningfulness is prior to truth'' \cite[p.~207]{fef08Axioms}.\footnote{
		The thought that determinateness is prior to truth was already expressed by Leitgeb in his semantic account of groundedness \cite{lei05Depends}, an axiomatization of which was later suggested by Schindler \cite{sch14Axioms}. However, the Leitgeb-Schindler approach imposes a drastic restriction, since groundedness \mrm G is `typed', in the sense that \mrm G provably does not apply to sentences containing \mrm G itself.
	}
A further incarnation of the Reinhardt-Feferman idea that determinateness can be defined in terms of truth has recently been developed by Castaldo and Nicolai \cite{cn25Classical}. They introduce the theory \ckf---the Classical Closure of \kfcs---in which the truth predicate plays both a logico-mathematical and a semantic role, and in which determinate sentences are those that are \textit{semantically} true or false.

In recent years, however, an alternative approach has emerged. Field \cite{fie22Power} and Fujimoto and Halbach \cite{fh24Classical} take admissibility---under the labels of \textit{strong classicality} and \textit{determinateness}, respectively---to be a primitive notion, not definable in terms of truth and falsity alone.\footnote{
		That a \textit{primitive} notion of determinateness is needed for Fujimoto and Halbach's purposes has been challenged by \cite{cn25Classical}, who argue that ``virtually all of the theoretical benefits sought for by Fujimoto and Halbach can be met by defining determinateness in terms of truth'' (p.~1043).
	}
Field's \INT extends the Strong Kleene system \pkf---Partial Kripke-Feferman truth---with a strong classicality predicate. Fujimoto and Halbach's \cd and \cdp add determinateness to a ``thoroughly classical'' theory of truth \cite[p.~220]{fh24Classical}. In both cases, the new predicate can be seen as embodying Kripke's notion of \textit{groundedness} from \cite{kri75Outline}.\footnote{
    Cf. \cite[p.~227]{fie22Power} and \cite[p.~222]{fh24Classical}.
}
And both theories show that axiomatizing admissibility together with truth yields systems which are conceptually rich and mathematically strong. Yet, the theories pull in opposite directions with respect to logic and truth. \INT is a \textit{nonclassical} theory of \textit{nonclassical} truth, i.e., it is formulated over a nonclassical logic and the logic governing the sentences it declares true is also nonclassical; \cd and \cdp are \textit{classical} theories of \textit{classical} truth. 

The current article continues this line of investigation, developing a third and alternative approach: it suggests axiomatizing determinateness over a \textit{classical} theory of \textit{nonclassical} truth, i.e., a theory that preserves classical logic, but treats truth itself in a nonclassical way. The theory will be labeled \kfdcs, for Kripke Feferman Truth with Consistency and Determinateness.

From a philosophical perspective, the ``most fundamental decision for the truth theorist is whether to sacrifice classical logic for transparency [of truth] or transparency for classical logic'' \cite[p.~259]{fh24Classical}.
We will not enter the discussion classical vs. nonclassical logic in this paper. Rather, we will sketch a comparison between \kfd and the classical alternative \cdp.\footnote{
    This will only be a sketch. A thorough philosophical analysis is left for a separate article. 
}
The theories \cdp and \kfd capture the same concept of determinateness, but they embody alternative concepts of truth. In effect, once classical logic is kept, a second crucial decision is forced. A consistent type-free truth predicate cannot be both complete, i.e. such that $\T\vphi\lor\T\neg\vphi$, and iterable, i.e. such that $\T\vphi\lra\T\T\vphi$. One of the two principles must be restricted and this is the crucial difference between \kfd and \cdp: the former retains iterability, the latter retains completeness.

\subsection{Iterability of \kfd vs.~Completeness of \cdp}\label{sse:CompVSIter}
	Fujimoto and Halbach suggest, as mentioned, a \textit{thoroughly classical} conception of truth, on which truth commutes with all connectives and quantifiers. Such a notion is fully characterized by the following axioms (see \cite[p.~220]{fh24Classical}):\footnote{
		The labels are those of \cite[p.~220]{fh24Classical}. We use a standard formal notation, explained in the next section.
	}
	\begin{enumerate}[label=\msf{T\arabic*}]
	\setcounter{enumi}{3} 
		\item $\forall \vphi\in\Sent \left(\T\neg \vphi\lra \neg \T\vphi\right)$
		\item $\forall \vphi,\psi\in\Sent \left(\T(\vphi\land \psi)\lra \T\vphi\land \T\psi\right)$
		\item $\forall\vphi\in\Fml^1 \left (\T(\forall v\vphi)\lra \forall t\,\T\vphi(t/v)\right)$.
	\end{enumerate} 
Iteration is then restricted to determinate sentences:
	\begin{enumerate}[label=\msf{T\arabic*}]
	\setcounter{enumi}{2} 
		\item $\forall t\bigl(\D \val t \ra (\T\ud\T t\lra\T\val t)\bigr)$.
	\end{enumerate}

This conception of truth is defended by pointing out that disquotational principles are not the only principles we rely on in our reasoning involving truth: compositional principles ``are just as important and used routinely without qualms'' \cite[p.~219]{fh24Classical}. Fujimoto and Halbach then provide three examples, similar in spirit, to show that we rely on tacit and often unnoticed uses of compositional principles. One of these examples involves the form of reasoning that Fujimoto \cite{fuj19Deflationism,fuj22Function} calls \textit{blind deductions}. For instance:
	\begin{quote}
		Alfred made a claim that was denied by Kurt (i.e., Kurt asserted
the negation of the claim); and everything Kurt claimed is true. Therefore not everything Alfred claimed is true. \hfill \cite[p.~220]{fh24Classical}
	\end{quote}
The validity of this reasoning relies on the consistency direction (left-to-right) of \msf{T4}: $x$ is not true if its negation is true. Note that we quantify over \textit{all} sentences, and we apply compositional principles `blindly', without knowing whether the sentences involved are determinate. 

The relevance of blind deductions can hardly be disputed. But they do not favor completeness over iterability.
On the one hand, blind deductions can just as well rely on the iteration of truth:
	\begin{quote}
        James claimed that some statements he read in an article on reflection principles are false. Everything James claims is true. Therefore, some statements in the article are false.
	\end{quote}
Besides the commutation of truth with connectives and quantifiers, the validity of this reasoning relies on the iteration principle that if `$x$ is false' is true, then $x$ is false, i.e., $\T\ud\F x\ra \F x$, where as usual we can define $\F x:\lra\T\ng x$. \cdp has this principle only for determinate sentences.

On the other hand, the blind deductions discussed in the literature (e.g.~\cite{fuj19Deflationism,fuj22Function,fh24Classical,cn25Classical}) involve only fairly uncontroversial instances of compositionality: it is clear that a conjunction is true precisely if both conjuncts are true, or that, \textit{pace} dialetheists, truth is consistent.
Yet, none of them relies on the \textit{completeness} of truth. And while a classical logician must accept that \vphi is \textit{true or not true}, they need not accept that \vphi is \textit{false} if it is not true. Indeed, if $\neg\T\vphi$ is read as a \textit{rejection} of the truth of \vphi, then there are cases discussed in the philosophical literature where the inference from $\neg\T\vphi$ to $\F\vphi$ is not warranted.\footnote{
    Incurvati and Schl\"oder \cite{is17Weak} distinguish between `strong' and `weak' rejection of a sentence \vphi. While the former amounts to asserting the falsity of $\vphi$, the latter does not. 
} 
We do not claim that these remarks are decisive---a thorough philosophical comparison is beyond the scope of this article---only that the case for completeness is no stronger than the case for iterability.

\subsection{Main results}

From a mathematical point of view, the results presented below confirm that a primitive notion of determinateness can significantly increase both the expressive power of a theory and its proof-theoretic strength. The key to the analysis is that \kfd contains, in effect, two truth predicates: `true' and `true and determinate'. The latter behaves like a Kripke-Feferman truth predicate, the former like a Kripke-Feferman truth predicate over a language containing the latter. In particular, in \kfd `true' is fully disquotational over `true and determinate'. This allows the interpretability of \kft (\cite{jkss99Proof, fuj11Autonomous}) via the mapping of the lower truth predicate to `true and determinate', and the upper truth predicate to `true' (§\ref{sec:LB}). The semantics via a double fixed-point construction (§\ref{sec:semantics}) tells the same story. Conversely, \kfdcs is interpreted in \kft extended with consistency and a principle of strong persistence (§\ref{sec:UB}), which turn out to be conservative (§\ref{sec:redKFTCSSP_to_KFTCS}). It then follows that
    \[
        \vert\kfd\vert = \vert\kfdcs\vert = 
        \vert\kft\vert = 
        \vert\idh2\vert = \vphi_{\vphi_{\varepsilon_0}0}0.
    \]
The analysis thereby locates the source of \kfd's strength compared to \cdp. The latter is mutually interpretable with \CT (Tarski's compositional truth theory) over \kfcs, and thus is arithmetically equivalent to the system \RA{<\varepsilon_{\varepsilon_0}} of ramified analysis up to any ordinal below $\varepsilon_{\varepsilon_0}$ \cite{fh24Classical}. Its proof-theoretic ordinal is thus $\vphi_{\varepsilon_{\varepsilon_0}}0$. The choice between completeness and iterability reappears one level up: just as \kf is significantly stronger than \CT, \kfd is significantly stronger than \cdp.

Finally, it is shown that the correspondence between \kfd and \kft is preserved when induction is restricted (§\ref{sec:add_res}). In particular, for internal induction we show\footnote{
    The theory \ids2 is obtained from \idh2 by restricting induction to formulae positive in the fixed-point operators (cf.~\cite{js96Some}).
}
    \[
        \vert\kfdint\vert = \vert\kfdintcs\vert = 
        \vert\kftint\vert = \vert\ids2\vert = \vphi_{\vphi_\omega0}0.
    \]
If internal induction is restricted to total formulae, then the resulting theories are such that 
    \[
        \vert\kfdtot\vert = \vert\kfdtotcs\vert = 
        \vert\kfttot\vert = \vert\idh1\vert = 
        \vert\kf\vert = \vphi_{\varepsilon_0}0.
    \]
If induction is restricted to the arithmetical vocabulary, the theory is model-theoretic conservative over Peano arithmetic.

\section{Preliminaries}
    Let $\lpa$ denote the language of Peano arithmetic (\PA) in some standard signature. Then
	\begin{itemize}
		\item $\Lmc_{\msf{KF_1}}:=\lpa\cup\{\Tz\}$,
		\item $\lkft:=\Lmc_{\msf{KF_1}}\cup\{\Tu\}$,
		\item $\lkfd:=\Lmc_{\msf{KF_1}}\cup\{\D\}$.
	\end{itemize}
\Lmc is used as metavariable for them. As usual, we define $\F_{i}x:=\T_{i}\ng x$ (for $i=0,1$). We often write \lkf instead of $\Lmc_{\msf{KF_1}}$, and for \lkf and \lkfd we often write \T instead of \Tz.
Formulae are generated by closing off under $\neg, \land, \forall$, with $\lor, \exists, \ra, \lra$ defined in the usual way according to classical logic. Greek letters $\vphi,\psi,\xi,\ldots$ range over formulae. The numeral corresponding to the number $n\in\omega$ is denoted by $\ovl n$. We fix a canonical G\"odel numbering and we take \PA as base theory, in which one can carry out a primitive recursive (p.r.) formalization of syntactic notions and operations. 

If $e$ is a term or a formula of $\Lmc$, the G\"odel number of $e$ is denoted by $\#e$ and $\cor{e}$ is the term representing $\#e$ in $\Lmc$. We take the following $\lpa$-predicates to abbreviate the equations for the characteristic functions of the sets terms, closed terms, variables, formulae, and sentences of $\Lmc$:
	\begin{IEEEeqnarray*}{L}
		\tm x \ (\ct x) := x \text{ is the gn of a (closed) term}; \\
		\vr x := x \text{ is the gn of a variable}; \\
		\Fml^n_{\Lmc}(x) \ (\Sent_\Lmc(x)) := x \text{ is the gn of an \Lmc-formula with at most } n \ (0) \text{ free distinct variables}
        \\
		\mrm{At}_{\Lmc}(x) := x \text{ is the gn of an atomic \Lmc-formula}.
	\end{IEEEeqnarray*}
For simplicity, we extend \Lmc with enough function symbols for elementary syntactic operations, including the usual ones on G\"odel numbers:
\begin{IEEEeqnarray*}{L+L+L}
    \ud= \text{ for } (\#t, \#s)\mapsto\#(t = s), & \ud R \text{ for } \#t\mapsto\#{R}{t} \text{ for } R \in\{\Tz,\Tu,\D\}
    \\
    \ng{} \text{ for } \#\vphi \mapsto\#\neg\vphi & \ud\land \text{ for } (\#\vphi,\#\psi)\mapsto\#(\vphi\land\psi), & \ud\forall \text{ for } (\#v,\#\vphi)\mapsto\#\forall v\vphi
\end{IEEEeqnarray*}
We also include a function symbol $\mrm{num}(x)$ for the standard numeral function, sending a number to the code of its numeral. Moreover, we employ a functional notation $x^\circ$ abbreviating the formula representing in $\lpa$ the evaluation function for closed terms of $\lpa$, such that ${\cor t}^\circ = t$ for closed terms $t$. For $x \in \Fml^1_\Lmc$, $x ( y/z )$, or just $x ( y )$, is a code of the sentence obtained by substituting the term $y$ for the free variable $z$ in the formula $x$. Moreover, for a formula $\vphi ( v )$, we define 
$\cor{\varphi(\dot{x})}:=\cor{\vphi(v)}( \num x / \cor v)$.

The expressions $\forall \vphi,\forall \psi,\forall \xi,\ldots$ abbreviate quantification over codes of formulae, and $\forall t,\forall s,\ldots$ abbreviate quantification over codes of closed terms. So for example, the expressions 
	\begin{IEEEeqnarray*}{LCL}
		\forall\vphi\in\Sent_\Lmc\left( \T_i\vphi \ra\neg\F_i\vphi\right)
        \\
        \forall \vphi\in\Fml^1_\Lmc \left(\T(\forall v\vphi)\lra \forall t\,\T\vphi(t)\right)
	\end{IEEEeqnarray*}
are shorthand for, respectively
	\begin{IEEEeqnarray*}{LCL}
		\forall x\left(\Sent_\Lmc(x) \ra (\T_i x \ra\neg\F_i x)\right)
		\\
        \forall x\forall y \left(\Fml^1_\Lmc(x)\land\mrm{Var}(y) \ra ( \T\all yx \lra \forall z (\ct z \ra \T x(z/y)))\right).
	\end{IEEEeqnarray*}

\section{The theories}

Theories in this paper will be formulated with `sentencehood axioms', i.e., with axioms saying that only \textit{sentences} can be true or determinate. The sole purpose of this technical assumption is to simplify some interpretability results. Moreover, we formulate the theories with the regularity axioms since we will also consider systems with restricted induction.

\begin{dfn}
	$\PA_\Lmc$ is \PA in the language \Lmc, including $\mrm{IND}_{\Lmc}$, where $\mrm{IND}_{\Lmc}$ is the full schema of induction for \Lmc:
	\begin{IEEEeqnarray*}{L+L}
		\hspace*{30mm} \bigl( \varphi ( 0 ) \wedge \forall x ( \varphi ( x ) \rightarrow \varphi ( x + 1 ) ) \bigr) \rightarrow \forall x \varphi ( x ), & \text{ for all \Lmc-formulae $\varphi ( x )$}.
	\end{IEEEeqnarray*}
\end{dfn}

\begin{dfn}
	$\msf{KF}_\mathit{i}$, for $\mathit{i}=1,2$, is the $\Lmc_{\msf{KF}_\mathit{i}}$-system expanding $\PA_{\Lmc_{\msf{KF}_\mathit{i}}}$ with the following axioms, where $k=0$ if $i=1$ and $k=0,1$ if $i=2$, and where \msf{KF_20} is an axiom only of \kft

\smallskip
\noindent
\textup{\bf KF-Axioms}
	\begin{enumerate}[label=$\msf{KF_\mathit{i}\arabic*}$]
	\setcounter{enumi}{-1}
        \item $\forall s\left(\left(\Tu\uTz s \lra \Tz\val{s}\right) \land \left(\Fu\uTz s \lra \neg\Tz\val{s}\right)\right)$.
		\item $\forall s\forall t\left((\T_k(\eqn st) \lra \val{s}=\val{t}) \land (\F_k(\eqn st) \lra \val{s}\neq\val{t})\right)$
		\item $\forall t \left((\T_k\val{t} \lra \T_k\ud\T_k t) \land (\F_k\val{t} \lra \F_k\ud\T_k t )\right)$
		\item $\forall \vphi\in\Sent_{\Lmc_{\kf_{k+1}}} \left(\T_k\vphi\lra\T_k\neg\neg\vphi\right)$
		\item $\forall \vphi,\psi\in\Sent_{\Lmc_{\kf_{k+1}}} \left((\T_k(\vphi\land \psi)\lra \T_k\vphi\land\T_k\psi) \land (\F_k(\vphi\land \psi)\lra \F_k\vphi \lor\F_k\psi)\right)$
		\item $\forall \vphi\in\Fml^1_{\Lmc_{\kf_{k+1}}} \left((\T_k(\forall v\vphi)\lra \forall t\,\T_k\vphi(t/v)) \land (\F_k(\forall v\vphi)\lra \exists t\,\F_k\vphi(t/v))\right)$
        \item $\forall \vphi\in\Fml^1_{\Lmc_{\kf_{k+1}}}\,\forall s\,\forall t\big(\val s= \val t \ra (\T_k\vphi(s) \lra \T_k\vphi(t))\big)$
        \item $\forall x\left(\T_k x\ra \Sent_{\lkft}(x)\right)$.
	\end{enumerate}
\end{dfn}

We mostly write \kf instead of $\msf{KF_1}$, and we write $\kf_\ltd$ for the system \kf in \ltd, where \ltd replaces $\Lmc_{\msf{KF}}$ in \msf{KF3}-\msf{KF6} and \msf{KF7} reads
\[
\forall x\left(\T x\ra\Sent_{\ltd}(x)\right).
\]

\begin{dfn}\label{dfn:kfd}
	\kfd is the \ltd-system extending $\kf_\ltd$ with the following axioms:
	
\smallskip
\noindent
\textup{\bf D-Axioms}
	\begin{enumerate}[label=$\msf{D\arabic*}$]
		\item $\forall s\forall t\left(\D(\eqn st)\right)$
		\item $\forall s \left(\D\val{s}\lra \D\ud\T s\right)$
		\item $\forall s \left(\D\val{s}\lra \D\ud\D s\right)$
		\item $\forall \vphi\in\Sent_{\ltd} \left(\D\vphi\lra \D\neg\vphi\right)$
		\item $\forall \vphi,\psi\in\Sent_{\ltd} \left(\D(\vphi\land \psi) \lra \left( (\D\vphi\land\D\psi) \lor (\D\vphi\land\F\vphi) \lor (\D\psi\land\F\psi) \right)\right)$
		\item $\forall\vphi\in\Fml^1_{\ltd} \left((\D(\forall v\vphi)\lra \forall t\,\D\vphi(t/v)) \lor (\exists t\,\D\vphi(t/v)\land\F\vphi(t/v))\right)$
		\item $\forall \vphi\in\Sent_{\ltd} \left(\D\vphi \ra (\F\vphi\lra\neg\T\vphi)\right)$
		\item $\forall s \left((\D\val{s}\lra \T\ud\D s) \land (\neg\D\val{s}\lra \F\ud\D s)\right)$.
        \item $(\forall \vphi(v)\forall s\forall t)\big(\val s= \val t \ra (\D\vphi(s) \lra \D\vphi(t))\big)$.
        \item $\forall x\left(\D x\ra \Sent_{\ltd}(x)\right)$.
	\end{enumerate}

\end{dfn}

\begin{dfn}
    For $S\in\{\kf_i,\kfd\}$, $S^-$ is defined as $S$ with induction restricted to $\lpa$.
\end{dfn}

\begin{dfn}
\hfill
    \begin{enumeratei}
    \item For $S\in\{\kfd, \kfd^-\}$, $S+\msf{Cons}$, or $S_\msf{cs}$ is the theory obtained by extending $S$ with 
	\begin{enumerate}[label=$\msf{Cons}$]
		\item \ $\forall \vphi\in\Sent_{\ltd}(\T\vphi\ra\neg\F\vphi)$.
	\end{enumerate}
    \item For $S\in\{\kf, \kf^-\}$, $S+\msf{Cons}$, or $S_\msf{cs}$ is the theory obtained by extending $S$ with \msf{Cons} by replacing $\Sent_{\ltd}$ with $\Sent_{\lkf}$.
    \item For $S\in\{\kft, \kft^-\}$, $S+\Const$ is the theory obtained by extending $S$ with 
	\begin{enumerate}[label=\Const]
		\item \ $\forall \vphi\in\Sent_{\lkf}(\Tz\vphi\ra\neg\Fz\vphi)\land\forall \vphi\in\Sent_{\ltt}(\Tu\vphi\ra\neg\Fu\vphi)$.
	\end{enumerate}
    \end{enumeratei}
\end{dfn}

It readily follows that the following compositional principles are derivable in \kfd:

\begin{lemma}\label{le:compositionalityKFD}
The following are derivable in \kfd:
    \begin{IEEEeqnarray}{L}
        \forall \vphi,\psi\in\Sent_{\ltd} \left((\T(\vphi\lor \psi)\lra \T\vphi\lor \T\psi) \land (\F(\vphi\lor \psi)\lra \F\vphi \land \F\psi)\right) \ztag{\msf{KF8}}
        \\
        \forall \vphi\in\Fml^1_{\ltd} \left((\T(\exists v\vphi)\lra \exists t\,\T\vphi(t/v)) \land (\F(\exists v\vphi)\lra \forall t\,\F\vphi(t/v))\right) \ztag{\msf{KF9}}
        \\
        \forall \vphi,\psi\in\Sent_{\ltd} \left(\D(\vphi\lor \psi) \lra \bigl((\D\vphi\land \D\psi) \land (\D\vphi\land \T\vphi) \lor (\D\psi \land \T\psi)\bigr)\right) \ztag{\msf{D11}}
        \\
        \forall\vphi\in\Fml^1_{\ltd} \left((\D(\exists v\vphi)\lra \forall t\,\D\vphi(t/v)) \lor (\exists t\,\D\vphi(t/v)\land\T\vphi(t/v))\right) \ztag{\msf{D12}}
    \end{IEEEeqnarray}
\end{lemma}

The crucial property for the proof-theoretic analysis to follow is that the truth predicate of \kfd is fully disquotational on $\T\land\D$:

\begin{lemma}\label{le:T-over-TcapD}
    $\kfd\vdash\forall t\big( (\T(\and{\ud\T t}{\ud\D t}) \lra\T\val t\land\D\val t) \land (\F(\and{\ud\T t}{\ud\D t}) \lra\neg\T\val t\lor\neg\D\val t) \big)$.
\end{lemma}

\begin{proof}
For the left conjunct:
	\begin{IEEEeqnarray*}{R}
		\T(\and{\ud\T t}{\ud\D t})
		\labeq{\msf{KF4}}\Lra 
		\T\ud\T t\land\T\ud\D t
		\labeq{\msf{KF2},\msf{D8}}\Lra
		\T\val t\land\D\val t.
	\end{IEEEeqnarray*}
For the right conjunct, first observe that 
    \[
        \F(\and{\ud\T t}{\ud\D t})\labeq{\msf{KF4}}\Lra \F\ud\T t\lor\F\ud\D t\labeq{\msf{KF2},\msf{D8}}\Lra\F\val t\lor\neg\D\val t.
    \]
Now we reason by cases. 
If $\D \val t$, then 
    \begin{IEEEeqnarray*}{R}
		\F\val t\lor\neg\D\val t\labeq{\msf{D7}, \D \val t}\Lra\neg\T\val t\lor\neg\D\val t
	\end{IEEEeqnarray*} 
If $\neg\D\val t$, then trivially $\F\val t\lor\neg\D\val t\Lra\neg\T\val t\lor\neg\D\val t$.
\end{proof}

\section{Semantics}\label{sec:semantics}
We provide fixed-point models for $\kfdrcs$, showing that Field's ``trick to run Kripke construction twice over'' \cite[p.~247]{fie22Power} can be used to obtain models for \kfdrcs and \kfdcs. 
We define a first jump for \T that does not contain clauses for sentences of the form $\neg\D s$.

\begin{dfn}\label{dfn_Tjump}
	Let $\mathscr{T}_i(x)$, for $1\leq i\leq 8$, be defined as follows:
	\begin{enumerate}[label={$\mathscr{T}_{\arabic*}(x) :=$}]\setlength{\itemindent}{8mm}
		\item $\exists s\exists t\left(\bigl(x=(\eqn st) \land \val s =\val t\bigr) \lor \bigl(x=(\ng{\eqn st}) \land \val s \neq\val t\bigr) \right)$
		\item $\Sent_{\ltd}(x)\land\exists s\left(\bigl(x=\ud\T s \land \T\val s\bigr) \lor \bigl(x=\ng{\ud\T s} \land \F\val s\bigr) \right)$
		\item $\Sent_{\ltd}(x)\land\exists s\left(x=\ud\D s \land \D\val s\right)$
		\item $\Sent_{\ltd}(x)\land\exists y\left(x = \ng{\ng}y \land \T y\right)$
		\item $\Sent_{\ltd}(x)\land\exists sr\left(x = \and sr \land \T s\land\T r\right)$
		\item $\Sent_{\ltd}(x)\land\exists sr\left(x = \ng{(\and sr)} \land (\F s\lor\F r)\right)$
		\item $\Sent_{\ltd}(x)\land\exists s\left(x = \all vs \land \forall y\T s(y)\right)$
		\item $\Sent_{\ltd}(x)\land\exists s\left(x = \ng{\all vs} \land \exists y\F s(y)\right)$.
	\end{enumerate}
Define $\mathscr{T}^-(x) := \bigvee_{1\leq j\leq 8} \mathscr{T}_j(x)$. Then, for $(\Nmc, X, Y)$ an \ltd structure where $\Nmc\vDash\PA$ and $X, Y$ interpret \D and \T, respectively, define:
	\begin{IEEEeqnarray*}{RCL}
		\Tmc^-_X(Y) &:=& \big\{a\in|\Nmc|:(\Nmc, X, Y)\vDash\mcr{T}^-(a) \big\}.
	\end{IEEEeqnarray*}
\end{dfn}

We next define an extension of the jump $\Tmc^-$, by adding a clause to the effect that sentences of the form $\neg\D t$ enter the jump whenever the value of $t$ is not determinate. 

\begin{dfn}
	Let 
	\begin{enumerate}
		\item[] $\mathscr{T}_9(x) := \Sent_{\ltd}(x)\land\exists s\left(x=\ng{\ud\D s} \land \neg\D\val s\right)$.
	\end{enumerate}
Define $\mathscr{T}(x):=\bigvee_{1\leq j\leq 9} \mathscr{T}_j(x)$ and $$\Tmc_X(Y) := \big\{a\in|\Nmc|:(\Nmc, X, Y)\vDash\mcr{T}(a) \big\}.$$ 
\end{dfn}

Finally, we define the jump for \D.

\begin{dfn}\label{dfn_Djump}
	Let $\mathscr{D}_i(x)$, for $1\leq i\leq 6$, be defined as follows:
	\begin{enumerate}[label={$\mathscr{D}_{\arabic*}(x) :=$}]\setlength{\itemindent}{8mm}
		\item $\exists s\exists t\left(x=(\eqn st)\right)$
		\item $\Sent_{\ltd}(x)\land\exists s\left(x=\ud\T s \land \D\val s \right)$
		\item $\Sent_{\ltd}(x)\land\exists s\left(x=\ud\D s \land \D\val s\right)$
		\item $\Sent_{\ltd}(x)\land\exists y\left(x = \ng y \land \D y\right)$
		\item $\Sent_{\ltd}(x)\land\exists sr\left(x = \and sr \land \left((\D s\land\D r) \lor (\D s\land\F s) \lor (\D r\land\F r)\right)\right)$
		\item $\Sent_{\ltd}(x)\land\exists s\left(x = \all vs \land \left(\forall y\D s(y) \lor \exists y(\D s(y)\land\F s(y))\right)\right)$.
	\end{enumerate}
Define $\mathscr{D}(x) := \bigvee_{1\leq j\leq 6} \mathscr{D}_j(x)$. Then define:
	\begin{IEEEeqnarray*}{RCL}
		\Dmc_Y(X) &:=& \big\{a\in|\Nmc|:(\Nmc, X, Y)\vDash\mcr{D}(a)\big\}.
	\end{IEEEeqnarray*}
\end{dfn}

Note that \D occurs only positively in $\Dmc_Y(X)$, and \T occurs only positively in both $\Tmc_X(Y)$ and $\Dmc_Y(X)$. Then, letting $(X,Y)\leq (X', Y'):\Lra X\subseteq X'\land Y\subseteq Y'$, the following monotonicity lemma can be proven in the usual way.

\begin{lemma}[Monotonicity]\label{le:monotonicity}
    The following jointly hold:
    \begin{enumerate}
        \item If $(X,Y)\leq (X', Y')$, then $\Dmc_Y(X)\subseteq\Dmc_{Y'}(X')$.
        \item If $(X,Y)\leq (X', Y')$, then $\Tmc^-_X(Y)\subseteq\Tmc^-_{X'}(Y')$.
        \item If $Y\subseteq Y'$, then $\Tmc_X(Y)\subseteq\Tmc_X(Y')$.
    \end{enumerate}
\end{lemma}

We now employ Field's `trick' of the double Kripkean fixed-point construction to obtain a model of \kfdcs.
We define the \textit{first run} as the sequence of interpretations for $D, T$ as follows,\footnote{
    Recall that the operator $\Tmc^-_X(Y)$ does not contain clauses for formulae of the form $\neg\D t$.
}
where for a set of sentences $Y$, we let $\ovl Y:=\{\vphi\mid\neg\vphi\in Y\}$:
    \begin{IEEEeqnarray*}{RCL+RCL}
        D_0 & := & \varnothing & T^-_0 & := & \varnothing
        \\
        D_{\alpha+1} & := & \Dmc_{T^-_\alpha}(D_\alpha) & T^-_{\alpha+1} & := & \Tmc^-_{D_\alpha}(T^-_\alpha)
        \\
        D_\lambda &:=& \bigcup_{\alpha<\lambda}D_\alpha & T^-_\lambda &:=& \bigcup_{\alpha<\lambda}T^-_\alpha.
    \end{IEEEeqnarray*}
It is easy to see that both sequences are weakly increasing and converging to a fixed point, let's call it $(D_\infty, T^-_\infty)$.
Since $\Tmc^-_X(Y)$ does not contain clauses for formulae of the form $\neg\D t$, the fixed-point $(D_\infty, T^-_\infty)$ will have `determinateness gaps', that is, some sentences of the form $\D t$ will be neither true nor false in the sense of $T^-_\infty$. In the \textit{second run}, these gaps are `closed off' by running $\Tmc_X(Y)$ over $(D_\infty, T^-_\infty)$. In particular, any determinateness ascription $\D t$ that was not declared true in the first run is declared false in the second.
More precisely, we reach a fixed-point $T_\infty$ via the sequence
    \begin{IEEEeqnarray*}{L+L+L}
        T_0 := T^-_\infty
        &
        T_{\alpha+1} := \Tmc_{D_\infty}(T_\alpha)
        &
        T_{\lambda} := \bigcup_{\alpha<\lambda}T_\alpha.
    \end{IEEEeqnarray*}
    
By definition, $T_\infty$ is such that $\neg\D t\in T_\infty$ iff $\D t\notin T^-_\infty$, and it can be shown that $(\Nbl, D_\infty, T_\infty)\vDash\kfdcs$.
We skip the proof since we'll shortly provide a more general characterization of the class of models of variants of \kfd.
Before that, though, let us observe that the fixed-point obtained above is identical to the one obtained by Field, in which only the truth jump is used and at every $\alpha$ determinateness is defined as truth-or-falsity.

\begin{lemma}\label{le:two_first_runs}
Let $(D'_\alpha,T'_\alpha)$ be defined like in the first run, but with $D'_{\alpha+1}:=T'_{\alpha+1}\cup\ovl{T'_{\alpha+1}}$, and let $(D'_\infty, T_\infty')$ be the fixed-point thus obtained. Then $(D'_\infty, T_\infty') = (D_\infty, T^-_\infty)$. In particular $D_\infty=T^-_\infty\cup\ovl{T^-_\infty}$.
\end{lemma}

\begin{proof}
    Both claims can be verified by a straightforward induction on stages. In particular, it is easy to see that $(D'_\alpha, T_\alpha') \subset (D_\infty, T^-_\infty)$ and, conversely, $(D_\alpha, T^-_\alpha) \leq (D'_\infty, T'_\infty)$.
\end{proof}

We now move on to the general characterization of models of variants of \kfd.

\begin{dfn}
    Given a model $\Nmc\vDash\PA$, let $\mrm{Sent}_{\ltd}^\Nmc:=\{a \mid \Nmc\vDash\mrm{Sent}_{\ltd}(a)\}$ denote the set of sentences in the sense of \Nmc. A set $X$ is called $\Nmc$-\emph{regular}, $X\in\mathrm{REG}_\Nmc$, iff $X\subseteq\mrm{Sent}^\Nmc_{\ltd}$ and $X$ ``is closed under substitution'', in the sense that the following holds:
	\begin{IEEEeqnarray*}{L}
		\text{if } \Nmc\vDash\mrm{Fml}^1_{\ltd}(\vphi)\land\mrm{ClTm}(a)\land\mrm{ClTm}(b)\land \mrm{val}(a) = \mrm{val}(b), \text{ then } \vphi(a)\in X \text{ iff } \vphi(b)\in X.
	\end{IEEEeqnarray*}
\end{dfn}

Just as for standard sentences, for $Y\in\mrm{REG}_\Nmc$, we let $\ovl Y:=\{a\mid\ng a\in Y\}$. Following \cite[Proposition 5.4]{can89Notes}, it can be shown that regularity is preserved under application of both \Tmc and \Dmc. 
Unless otherwise specified, from now on we restrict ourselves to regular sets.

\begin{prop}\label{pr:modelcharacterization}
Let $\Nmc\vDash\PA$. Then:
\begin{enumeratei}
    \item $(\Nmc, X, Y)\vDash\kfdr$ iff the following jointly hold:
    \begin{enumerate}\setlength{\itemindent}{3mm}
        \item $X\subseteq Y\cup\ovl Y$ and $X\cap(Y\cap\ovl Y) = \emptyset$,
        \item $X=\Dmc_Y(X)$,
        \item $Y = \Tmc_X(Y)$.
    \end{enumerate}
    \item $(\Nmc, X, Y)\vDash\kfdrcs$ iff (a)-(c) hold along with
    \begin{enumerate}\setlength{\itemindent}{3mm}
        \item[(d)] $Y\cap\ovl Y = \emptyset$.
    \end{enumerate}
    \item If $\Nmc = \Nbl$, we can replace \kfdr with \kfd in (i) and \kfdrcs with \kfdcs in (ii).
\end{enumeratei}
\end{prop}

\begin{proof}
    Let $(\Nmc, X, Y)$ be a model of $\kfdr$. $(a)$ follows by \msf{D7}. The proofs for $(b)$ and $(c)$ are routine, reasoning by cases on $a\in\Sent^\Nmc$. For example, if $\Nmc\vDash a = \vphi\land\psi$, then we have:
	\begin{IEEEeqnarray*}{L}
		\vphi\land\psi \in X 
		\labeq{\msf{D4,D5}}\Lra (\vphi,\psi\in X) \lor (\neg \vphi\in X\cap Y) \lor (\neg \psi\in X\cap Y)
		\labeq{\mcr{D}_4,\mcr{D}_5}\Lra \vphi\land\psi \in \Dmc_{Y}(X). 
	\end{IEEEeqnarray*}   

As for the converse direction, it is straightforward to check that all axioms of $\kfdr (+\,\Cons)$ are satisfied in $(\Nmc, X, Y)$ if $(a)$-$(c)$ (and $(d)$) are satisfied. We only show a few examples, beginning with \msf{D5}. We have $\vphi\land\psi\in X=\Dmc_Y(X)$ iff $\vphi, \psi\in X$ or $\neg\vphi\in X\cap Y$ or $\neg\psi\in X\cap Y$. Consider \msf{D7} and let $\vphi\in X$. The conclusion follows immediately by $(a)$. As for \msf{D8}, we have $s\in X \, (s\notin X)$ iff $\D s\in\Tmc_X(Y)=Y \, (\neg\D s\in\Tmc_X(Y)=Y)$. 
Claims (ii) and (iii) are obvious.
\end{proof}

The double fixed-point construction introduced above can be run over arbitrary models of \PA, and the following lemma is clear.

\begin{lemma}\label{le:consistency}
Let $X$ be any set of sentences in the sense of \Nmc. If $Z$ is consistent and $Z\subseteq\Tmc_X(Z)$, then $\Tmc_X(Z)$ is consistent. The same holds with respect to $\Tmc^-_X$. 
Hence, the least fixed points of $\Tmc_X$ and $\Tmc^-_X$ are consistent, and so are $T^-_\infty$ and $T_\infty$.
\end{lemma}

We next show that models obtained by the double construction satisfy (a)-(d) of \Cref{pr:modelcharacterization}. This establishes that $(\Nbl,D_\infty,T_\infty)\vDash\kfdcs$, and that $\kfdrcs$ is model-theoretically conservative over \PA (\Cref{thm:FPmodelKFDcs}). We now fix an arbitrary \Nmc and dependence upon it will be tacit from context.

\begin{dfn}
Let $\Kmc(X, Y) := (\Dmc_Y(X), \Tmc^-_X(Y))$, and $\UPm(X, Y):=\sup\{(X_\alpha, Y_\alpha)\mid\alpha\in ON\}$, where 
		\begin{IEEEeqnarray*}{L+L+L}
			(X_0, Y_0)=(X, Y), & (X_{\alpha+1},Y_{\alpha+1})=\Kmc(X_\alpha, Y_\alpha), & (X_\lambda, Y_\lambda)=\sup\{(X_\beta, Y_\beta)\mid\beta<\lambda\}.
		\end{IEEEeqnarray*}
We analogously define $\UP_X(Y)$ as the upward closure of $\Tmc_X(Y)$.
\end{dfn}

Note that \Cref{le:monotonicity} entails that \Kmc is monotone. Then the following readily follows:

\begin{prop}\label{pr:upwardclosure}
    \hfill
    \begin{enumeratei}
        \item If $(X, Y)\leq\Kmc(X, Y)$, then $\UPm(X, Y)$ is the least fixed-point of $\Kmc$ above $(X, Y)$.
        \item If $Y\subseteq\Tmc_X(Y)$, then $\UP_X(Y)$ is the least fixed-point of $\Tmc_X$ above $Y$.
    \end{enumeratei}
\end{prop}

\begin{obs}\label{obs:soundnessfirstrun}
    Let $(X, Y)=\Kmc(X, Y)$. Then $Y\subseteq\Tmc_X(Y)$.
\end{obs}

\begin{proof}
    Immediate, since $Y=\Tmc^-_X(Y)\subseteq\Tmc_X(Y)$.
\end{proof}

The next lemma states in general terms that the intersection of true-and-determinate sentences is invariant across the two runs as long as determinateness does not contain gluts.

\begin{lemma}\label{le:intersectionsonD}
    Let $(X, Y)=\Kmc(X, Y)$ be such that $X\cap(Y\cap \ovl Y)=\emptyset$. Let $Z=\UP_X(Y)$, and write $Z_\alpha$ instead of $Y_\alpha$. Then $X\cap Y = X\cap Z$.
\end{lemma}

\begin{proof}
    It suffices to show that $X\cap Z_\alpha\subseteq Y$ for all $\alpha$ by transfinite induction on $\alpha$. Since the claim is obvious for $\alpha=0$ or limit, let $\alpha=\beta+1$ and suppose $\vphi\in X\cap Z_{\beta+1}$. The desired conclusion can be obtained reasoning by cases depending on the shape of \vphi. For example, if $\vphi=\T s$, then $\val s\in X\cap Z_\beta$. By IH, $\val s\in X\cap Y$ and thus $\T s\in X\cap Y$. Note that the case where $\vphi=\neg\D s$ cannot arise: from $\neg\D s\in Z_{\beta+1}$ we derive that $\val s\notin X$, and hence $\{\D s, \neg\D s\}\cap X=\emptyset$. If $\vphi=\psi\land\xi$, then $\vphi\in X\cap Z_{\beta+1}$ entails that $\psi,\xi\in Z_\beta$; moreover, from $\vphi\in X$ it follows that either $\psi,\xi\in X$, hence $\psi,\xi\in Y$ by IH and thus $\vphi\in Y$, or one conjunct, say $\psi$, is such that $\psi\in X$ and $\neg\psi\in Y$. But then $\psi\in Y$ by IH, contrary to the assumption.
\end{proof}

\begin{lemma}\label{le:characterization}
    Let $(X, Y)=\Kmc(X, Y)$ be such that $X\subseteq Y\cup \ovl Y$ and $X\cap(Y\cap \ovl Y)=\emptyset$, and let $Z=\UP_X(Y)$. Then 
    \begin{enumeratei}
        \item $(X, Z)$ satisfies (a)-(c) of \Cref{pr:modelcharacterization}. Hence, $(\Nmc, X, Z)\vDash\kfdr$.
        \item $(\Nmc, X, Z)\vDash\kfdrcs$ iff, in addition, $Y\cap\ovl Y=\emptyset$.
        \item If $\Nmc=\Nbl$ is standard, (i) and (ii) hold for \kfd and \kfdcs, respectively.
    \end{enumeratei}
\end{lemma}

\begin{proof}
(i) Condition (a) -- $X\subseteq Z\cup \ovl Z$ and $X\cap(Z\cap \ovl Z)=\emptyset$ -- follows from the hypothesis via \Cref{le:intersectionsonD}. As for (c) -- $Z=\Tmc_X(Z)$ --, by \Cref{obs:soundnessfirstrun} we have $Y\subseteq\Tmc_X(Y)$; thus, by \Cref{pr:upwardclosure}, $Z$ is the least fixed-point of $\Tmc_X$ above $Y$. As for (b) -- $\Dmc_Z(X)=X$ --, by \Cref{le:intersectionsonD} we have that $A:=\{\vphi\in X\mid \neg\vphi\in Z\} = \{\vphi\in X\mid \neg\vphi\in Y\}=:B$. Since $\Dmc_Z(X)$ depends on $Z$ only through $A$, we have $\Dmc_Z(X)=\Dmc_Y(X)=X$.

(ii) Given (i), it suffices to show that $Z\cap \ovl Z=\emptyset$ iff $Y\cap \ovl Y=\emptyset$. The left-to-right direction is obvious, since $Y\subseteq Z$. For the converse direction, it suffices to observe that by \Cref{le:consistency} consistency is preserved under $\Tmc_X$ for any $Y$ such that $Y\subseteq \Tmc_X(Y)$. (iii) is obvious.
\end{proof}

Since the double construction over $(\emptyset, \emptyset)$ yields a structure satisfying the conditions of \Cref{le:characterization}, we obtain the following:

\begin{theorem}\label{thm:FPmodelKFDcs}
    (i) $(\Nbl, D_\infty, T_\infty)\vDash\kfdcs$. (ii) $\kfdrcs$ is model-theoretic conservative over \PA.
\end{theorem}

Note that the the double construction starting from the empty set does not yield the least fixed-point model of \kfdcs, for $\Tmc_X(Y)$ is antitone in $X$ and treats $X$ classically. For example, given a truth-teller $\tau$, we have $\neg\D\tau\in T_\infty$, but $\neg\D\tau$ is not true is not true in any fixed-point satisfying $\D\tau$. Nonetheless, it can be shown that $(\Nbl, D_\infty, T_\infty)$ has certain minimality properties. We write $(D_\alpha,T^-_\alpha)$ for the stages of the first run on $(\emptyset, \emptyset)$ and $T_\alpha$ for the stages of the second run.

\begin{lemma}\label{le:closed_pairs}
If $\Dmc_Y(X)\subseteq X$ and $\Tmc_X(Y)\subseteq Y$, then $D_\infty\subseteq X$ and $T^-_\infty\subseteq Y$.
\end{lemma}

\begin{proof}
By induction on $\alpha$, we show that $D_\alpha\subseteq X$ and $T^-_\alpha\subseteq Y$. The cases $\alpha=0$ and $\alpha$ limit are trivial. For the successor step:
 \begin{IEEEeqnarray*}{L+L+L+L}
    && D_{\alpha+1}=\Dmc_{T^-_\alpha}(D_\alpha)\labeq{\text{IH, \Cref{le:monotonicity}}}\subseteq\Dmc_Y(X)\subseteq X,
    \\
    && T^-_{\alpha+1}=\Tmc^-_{D_\alpha}(T^-_\alpha)\labeq{\text{IH, \Cref{le:monotonicity}}}\subseteq\Tmc^-_X(Y)\subseteq\Tmc_X(Y)\subseteq Y. & \qedhere
 \end{IEEEeqnarray*}
\end{proof}

\begin{lemma}\label{le:consistencyfirstrun}
    Let $X\supseteq D_\infty$ and $Z=\Tmc_X(Z)$. Then the following hold.
    \begin{enumeratei}
        \item $T^-_\infty\subseteq Z$ and $\ovl{T^-_\infty}\cap Z=\emptyset$.
        \item $D_\infty\cap Z=D_\infty\cap T_\infty^-$. In particular, $D_\infty\cap T^-_\infty = D_\infty\cap T_\infty$.
    \end{enumeratei}
\end{lemma}

\begin{proof}
(i) For $T^-_\infty\subseteq Z$, we show $T^-_\alpha\subseteq Z$ by induction on $\alpha$, as in the proof of \Cref{le:closed_pairs}: since $D_\alpha\subseteq D_\infty\subseteq X$, IH and \Cref{le:monotonicity} give $T^-_{\alpha+1}=\Tmc^-_{D_\alpha}(T^-_\alpha)\subseteq\Tmc^-_X(Z)\subseteq\Tmc_X(Z)=Z$. For the claim $\ovl{T^-_\infty}\cap Z=\emptyset$, by a straightforward induction on $\alpha$ one first shows that $\vphi\in T^-_\alpha$ entails $\neg\vphi\notin Z$, and then obtains the desired conclusion by closure under double negation of $Z$.

(ii) Let $\vphi\in D_\infty\cap Z$, which by the equivalence $D_\infty = T^-_\infty\cup\ovl{T^-_\infty}$ (\Cref{le:two_first_runs}) yields $\vphi\in (T^-_\infty\cup\ovl{T^-_\infty})\cap Z$. By (i), $\ovl{T^-_\infty} \cap Z = \emptyset$, hence we conclude $\vphi\in T^-_\infty$.
\end{proof}

By \Cref{pr:modelcharacterization}, models of \kfd satisfy the hypothesis of \Cref{le:closed_pairs}. We thus obtain:

\begin{crl}[Minimality of $D_\infty$]
     $D_\infty\subseteq X$ for every model $(\mathbb N,X,Y)$ of \kfd.
\end{crl}

Finally, we want to show that $(D_\infty, T_\infty)$ is the quasi-least model of \kfdcs, in the sense that it is the unique $(X, Z)$ such that $X$ is the least fixed point of $\Dmc_{Z}$ and $Z$ is the least fixed point of $\Tmc_{X}$.

\begin{lemma}\label{le:quasi-minimality}
\hfill
\begin{enumeratei}
  \item $T_\infty$ is the least fixed point of $\Tmc_{D_\infty}$.
  \item $D_\infty$ is the least fixed point of $\Dmc_Y$, for every $(X, Y)$ such that $X\supseteq D_\infty$ and $Y=\Tmc_X(Y)$.
\end{enumeratei}
\end{lemma}

\begin{proof}
Claim (i) follows from \Cref{pr:upwardclosure} and \Cref{le:consistencyfirstrun}(i) with $X=D_\infty$. For (ii), let $X\supseteq D_\infty$ and $Y=\Tmc_X(Y)$. Let $L$ be the least fixed point of $\Dmc_Y$. We need to show that $D_\infty = L$. Let us first show that $D_\infty$ is a fixed point of $\Dmc_Y$. Given that $Y$ and $T^-_\infty$ agree on which sentences in $D_\infty$ are false by \Cref{le:consistencyfirstrun}, and given that
$\Dmc_Y(X)$ depends on $Y$ only through $\{\psi\in X:\neg\psi\in Y\}$, we get $\Dmc_Y(D_\infty)=\Dmc_{T^-_\infty}(D_\infty)=D_\infty$. So $D_\infty$ is a fixed point of $\Dmc_Y$. To show that $D_\infty\subseteq L$, we show that, for all $\alpha$, $D_\alpha\subseteq L$. We let $\alpha=\beta+1$ and reason by cases depending on the shape of \vphi. If \vphi is a truth atom $\T t$, then $\val t\in D_\beta$. By IH, $\val t\in L$ and therefore $\T t\in L$ by fixed-point properties of $L$. The argument is similar if \vphi is a positive determinateness atom $\D t$ or a negation $\neg\psi$. If $\vphi=\psi\land\xi$, then either $\psi,\xi\in D_\beta$ and we are done by IH, or one conjunct, say $\psi$, is false and determinate, i.e., $\neg\psi\in D_\beta\cap T_\beta^-$. By Lemma \ref{le:consistencyfirstrun} and IH, this yields $\neg\psi\in L\cap Y$, hence $\vphi\in L$.
\end{proof}

\begin{prop}[Quasi-minimality]\label{pr:quasi-minimality}
    $(D_\infty,T_\infty)$ is the unique pair $(X,Y)$ of sets of sentences such that $X$ is the least fixed point of $\Dmc_Y$ and $Y$ is the least fixed point of $\Tmc_X$.
\end{prop}

\begin{proof}
By \Cref{le:quasi-minimality}, $(D_\infty,T_\infty)$ has the property. For uniqueness, suppose $(X,Y)$ is such that $X$ is the least fixed point of $\Dmc_Y$ and $Y$ is the least fixed point of $\Tmc_X$. In particular, $\Dmc_Y(X)=X$ and $\Tmc_X(Y)=Y$, so $D_\infty\subseteq X$ by \Cref{le:closed_pairs}. Since $X\supseteq D_\infty$ and $Y=\Tmc_X(Y)$, \Cref{le:quasi-minimality}(ii) yields that $D_\infty$ is the least fixed point of $\Dmc_Y$, that is, $X=D_\infty$. Hence $Y$ is the least fixed point of $\Tmc_{D_\infty}$, hence $Y=T_\infty$ by \Cref{le:quasi-minimality}(i).
\end{proof}

\section{Reduction of \kft to \kfd}\label{sec:LB}
	In this section, we show that \kfd proves every arithmetical theorem of theory $\idh2$ of twice iterated fixed-points. The systems $\widehat{\msf{ID}}_n$ of finitely many times iterated fixed-points were introduced by Feferman in \cite{fef82Iterated}, where it was established that $\widehat{\msf{ID}}_n$ has proof-theoretic ordinal $\alpha_n$, for $\alpha_n$ the $n+1$-th element of the canonical fundamental sequence for $\Gamma_0$ with $\alpha_0:=\varepsilon_0$ and $\alpha_{n+1}:=\vphi_{\alpha_n}(0)$. Later, a transfinite progression of systems of Kripke-Feferman truth was introduced and studied by \cite{jkss99Proof}---see \cite{fuj11Autonomous} for further analysis and \cite[Ch. 8]{can96Logical} for a related approach.\footnote{
    In \cite{jkss99Proof} the systems are called $\msf{SRT}_\alpha$, for `self-reflective truth', but we use the label $\msf{KF}_\alpha$ from \cite{fuj11Autonomous}.
}
They show that $\msf{KF}_\alpha$ has the same proof-theoretic strength as the theory $\widehat{\msf{ID}}_\alpha$ of $\alpha$-many times iterated fixed-points.\footnote{
    Specifically, \cite{jkss99Proof} show that
    \[
        \vert\msf{KF}_\alpha\vert = \vert\widehat{\msf{ID}}_\alpha\vert = \vphi1\gamma_n(\vphi1\gamma_{n-1}(\ldots\vphi1\gamma_1(\alpha|m)) \ldots ))
    \]
where: (i) $\vphi \alpha\beta$ is the ternary Veblen function, (ii) $\varepsilon(\delta)$ denotes the least epsilon number above $\delta$, (iii) $(\delta|m)$ is inductively defined by $(\delta|0)=\varepsilon(\delta)$ and $(\delta|k+1) = \vphi_{\delta|k}0$, and (iv) $\alpha = \omega^{1+\gamma_n} + \omega^{1+\gamma_{n-1}} + \ldots + \omega^{1+\gamma_{1}} + m$ with $\gamma_n\geq\ldots\geq\gamma_1$.
}
We can thus obtain our result by constructing a truth-definition, in the sense of \cite{fuj10Relative}, of the system \kft in \kfd.

We begin by defining a translation of \ltt into \ltd.

\begin{dfn}\label{dfn_flat}
By the recursion theorem, let $k$ be a function satisfying the following conditions:
	\begin{IEEEeqnarray}{L}
		k x:=
		\begin{cases}
			x & \text{if $x\in\mrm{At}_{\lpa}$}
			\\
			\and{\ud\T\ud k y}{\ud\D\ud k y} & \text{if $x = \uTz y$, for $y\in\mrm{Term}$}
			\\
			\ud\T\ud k y & \text{if $x = \uTu y$, for $y\in\mrm{Term}$}
			\\
			\ng k y & \text{if $x=\ng y\in\mrm{Fml}_{\ltt}$}
			\\
			\and{ky}{kz} & \text{if $x=\and yz\in\mrm{Fml}_{\ltt}$}
			\\
			\all{v}{ky} &  \text{if $x=\all vy\in\mrm{Fml}_{\ltt}$}
			\\
			0 & \text{otherwise}
		\end{cases}
		\label{eq_k}
	\end{IEEEeqnarray}
As usual, $\ud k$ is an arithmetical representation of the syntactic operation of applying $k$ to terms. We require that, for $s\in\ClTm$,
	\begin{IEEEeqnarray}{L}
		\PA\vdash\ClTm(\ud ks) \land \val{(\ud ks)} = k(\val s). \label{eq_PAk}
	\end{IEEEeqnarray}
Let $\flat$ be a translation of $\ltt$ in \ltd such that $\flat$ preserves the arithmetical vocabulary, commutes with logical operations and quantifiers, and is such that
	\begin{IEEEeqnarray*}{L+C+L}
		\Tz^\flat x := \T kx \land \D kx
		& \text{and} &
		\Tu^\flat x := \T k x.
	\end{IEEEeqnarray*}
\end{dfn}

By convention, $\ng x=0$ for any $x\notin\Fml_\ltt$; hence $k\ng x=\ng k x$ for all $x$. This allows the following simplification relative to the translations of falsity ascriptions:
	\begin{IEEEeqnarray}{L}
		\Fz^\flat x \Lra \T k\ng{x}\land \D k\ng{x} \Lra \F kx\land\D \ng kx \labeq{\msf{D4}}\Lra \F kx\land\D kx; \label{eq_Faux}
		\\
		\Fu^\flat x \Lra \T k\ng{x} \Lra \F kx. \label{eq_F2aux}
	\end{IEEEeqnarray}

\begin{lemma}\label{le_kf_21-5inkfd}
	The $\flat$-translations of \msf{KF_21}-\msf{KF_27} are derivable in $\kfd$.
\end{lemma}

\begin{proof}
We begin with \msf{KF_27}. Its $\flat$-image is equivalent to $\forall x\left(\T kx \ra \Sent_{\lkft}(x)\right)$. By \msf{KF7}, $\forall x(\T kx \ra \Sent_{\lkfd}(kx))$. Provably in \PA, $\Sent_{\lkfd}(kx)$ only if $x\in\Sent_\lkft$. Thus, $\forall x(\T kx \ra \Sent_{\lkft}(x))$.

As for \msf{KF_21}-\msf{KF_26}, the argument for the case $\T_i=\Tu$ is straightforward, given that $\Tu$ is mapped to \T and \kfd has the corresponding axioms. For example, the $\flat$-translation of \msf{KF_22} can be shown to be derivable as follows:
	\begin{IEEEeqnarray*}{C}
		\Tu^\flat \val t 
		\Lra
		\T k\val t 
		\labeq{\eqref{eq_PAk}}\Lra 
		\T \val{(\ud kt)} 
		\labeq{\msf{KF2}}\Lra 
		\T\ud\T\ud kt 
		\Lra
		\T k\ud{\Tu} t 
		\Lra
		\Tu^\flat\ud{\Tu} t;
		\\
		[2mm]
		\Fu^\flat \val t 
		\labeq{\eqref{eq_F2aux}}\Lra
		\F k\val t 
		\labeq{\eqref{eq_PAk}}\Lra 
		\F \val{(\ud kt)}
		\labeq{\msf{KF2}}\Lra 
		\F \ud\T\ud kt 
		\Lra
		\F k\uTu t 
		\labeq{\eqref{eq_F2aux}}\Lra
		\Fu^\flat\uTu t.
	\end{IEEEeqnarray*}
Other cases are treated similarly.

As for $\T_i=\Tz$, axioms \msf{KF_21} and \msf{KF_23} are obvious, so we begin with \msf{KF_22}. Observe that
	\begin{IEEEeqnarray}{L}
		\D k\val t 
		\labeq{\eqref{eq_PAk}}\Lra 
		\D\val{(\ud kt)} 
		\labeq{\msf{D2, D3}}\Lra
		\D\ud\T\ud kt \land\D\ud\D\ud kt 
		\labeq{\msf{D2, D3, D5}}\Lra
		\D(\and{\ud\T\ud kt}{\ud\D\ud kt}). \label{eq_Daux}
	\end{IEEEeqnarray}
We then reason as follows:
	\begin{IEEEeqnarray*}{R}
		\Tz^\flat \val t 
		\labeq{\eqref{eq_k},\eqref{eq_PAk}}\Lra
		\T \val{(\ud kt)}\land \D\val{(\ud kt)}
		\labeq{\msf{KF2, D8}}\Lra 
		\T\ud\T\ud kt\land\T\ud\D\ud kt 
		\labeq{\msf{KF4,D8,\eqref{eq_Daux}}}\Lra
		\T (\and{\ud\T\ud kt}{\ud\D\ud kt})\land\D (\and{\ud\T\ud kt}{\ud\D\ud kt})	
		\\
		\labeq{\eqref{eq_k}}\Lra
		\T k\uTz t\land\D k\uTz t 
		\Lra 
		\Tz^\flat\uTz t.
		\\[2mm]
		\Fz^\flat \val t 
		\labeq{\eqref{eq_k},\eqref{eq_PAk},\eqref{eq_Faux}}\Lra
		\F\val{(\ud kt)}\land \D k{\val t}
		\labeq{\msf{KF2}}\Lra 
		\F\ud\T\ud kt \land \D k{\val t}
		\labeq{\msf{KF4,\eqref{eq_Daux}}}\Lra
		\F (\and{\ud\T\ud kt}{\ud\D\ud kt}) \land \D (\and{\ud\T\ud kt}{\ud\D\ud kt})
		\\
		\labeq{\eqref{eq_k}}\Lra
		\F k\uTz t \land \D k\uTz t 
		\labeq{\eqref{eq_Faux}}\Lra 
		\Fz^\flat\uTz t.
	\end{IEEEeqnarray*}

We verify the falsity conjunct of \msf{KF_24}---the other is similar. Fix $\vphi,\psi\in\Sent_{\lt}$, and observe
	\begin{IEEEeqnarray*}{RCL}
		\Fz^\flat(\vphi\land\psi) &\labeq{\eqref{eq_Faux}}\Lra& \F (k\vphi\land k\psi)\land\D(k\vphi\land k\psi) \labeq{\msf{KF4}} \Lra (\F k\vphi\lor\F k\psi)\land\D(k\vphi\land k\psi);
		\\[2mm]
		\Fz^\flat\vphi\lor\Fz^\flat\psi &\labeq{\eqref{eq_Faux}}\Lra& (\F k\vphi\land\D k\vphi)\lor(\F k\psi\land\D k\psi).
	\end{IEEEeqnarray*}
The implication $\Fz^\flat\vphi\lor\Fz^\flat\psi \ra \Fz^\flat(\vphi\land\psi)$ is then immediate by logic and \msf{D5}. For the converse, $\D(k\vphi\land k\psi)$ entails $\D k\vphi\lor\D k\psi$. Then, the claim is obvious if $\D k\vphi\land\F k\vphi$ or $\D k\psi\land\F k\psi$. If $\D k\vphi\land\F k\psi$, then from $\D(k\vphi\land k\psi)$ by \msf{D5} we derive $\D k\psi\lor\F k\vphi$. Similarly if $\D k\psi\land\F k\vphi$.

The argument for \msf{KF_25} is similar.
\end{proof}

\begin{lemma}\label{le_kf20inkfd}
	The $\flat$-translation of $\msf{KF_20}$ is derivable in $\kfd$.
\end{lemma}

\begin{proof}
	This is a key axiom, for which \Cref{le:T-over-TcapD} is crucial. For the left conjunct:
	\begin{IEEEeqnarray*}{R}
		\Tu^\flat \uTz t
		\Lra 
		\T k\uTz t
		\Lra
		\T (\and{\ud\T\ud kt}{\ud\D\ud kt})
		\labeq{\Cref{le:T-over-TcapD},\eqref{eq_PAk}}\Lra
		\T k \val t\land\D k \val t \Lra
		\Tz^\flat \val t
	\end{IEEEeqnarray*}
For the right conjunct:
	\begin{IEEEeqnarray*}{L+L+L+L}
		&& \Fu^\flat\uTz t
		\labeq{\eqref{eq_F2aux}}\Lra 
		\F k{\uTz}t
		\Lra 
		\F (\and{\ud\T\ud kt}{\ud\D\ud kt})
        \labeq{\Cref{le:T-over-TcapD},\eqref{eq_PAk}}\Lra
        \neg\T k\val t\lor\neg\D k\val t \Lra \neg\Tz^\flat\val t &  \qedhere
	\end{IEEEeqnarray*}
\end{proof}

\begin{theorem}\label{th_LB}
	The translation $\flat$ is a truth-definition of \kft in \kfd. Hence, $\vert\kft\vert=\alpha_2\leq\vert\kfd\vert$.
\end{theorem}

\section{Reduction of \kfdcs to \kftcssp}\label{sec:UB}

We prove that the lower bound is sharp, and that consistency is \lpa-conservative over \kfd. To this end, we begin by providing an \lpa-interpretation of \kfdcs in an extension of \kftcs. Specifically, \kftcs will be expanded with a principle of \textit{strong persistence}. In the next section, we then establish that this extension of \kftcs is interpretable in \kftcs. Let us first recall that, by a theorem of \cite{fuj11Autonomous}, consistency is conservative over iterations of \kf.

\begin{theorem}[{\cite[Theorem 16]{fuj11Autonomous}}]\label{th_FujCon}
	\kfp\alpha and $\kfp\alpha+\msf{Cons_\alpha}$ prove the same arithmetical sentences, where \msf{Cons_\alpha} is an axiom postulating the consistency of every truth predicate of \kfp\alpha.\footnote{
    Fujimoto's techniques for proving the conservativity of \msf{Cons_\alpha} will be used in §\ref{sec:add_res} to obtain the conservativity of consistency over \kft with internal induction.
}
\end{theorem}

We begin by defining a translation of \ltd in \ltt. Coded occurrences of \T will not be treated uniformly: occurrences within the scope of \D will be mapped to \Tz, while the remaining ones will be mapped to \Tu. 

\begin{dfn}\label{dfn_natural}

Let $g$ and $f$ be p.r.~functions satisfying the following conditions:
	\begin{IEEEeqnarray}{L}
		gx:=
		\left\{\def\arraystretch{1.2}
		\begin{array}{@{}l@{\hspace{10mm}}l@{\hspace{10mm}}l@{}}
			x & \text{if $x\in\mrm{At}_{\lpa}$} & x
			\\
			\ud{\Tz}\ud gy & \text{if $x = \ud\T y$ for $y\in\Term$} & \ud{\Tu}\ud fy
			\\
			\dor{\uTz\ud gy}{\uFz\ud gy} & \text{if $x = \ud\D y$ for $y\in\Term$} & \dor{\uTz\ud gy}{\uFz\ud gy}
			\\
			\ng gy & \text{if } x=\ng y\in\mrm{Fml}_{\lkfd} & \ng fy
			\\
			\and{gy}{gz} & \text{if $x=\and yz\in\mrm{Fml}_{\lkfd}$} & \and{fy}{fz}
			\\
			\all v{gy} & \text{if $x=\all vy\in\mrm{Fml}_{\lkfd}$} & \all v{fy}
			\\
			0 & \text{otherwise} & 0
		\end{array}
		\right\}
		=:fx
		\label{eq_gf}
	\end{IEEEeqnarray}
As above, $\ud g$ and $\ud f$ are arithmetical representations of the syntactic operation of applying $g$ and $f$, respectively, to terms. We have, for $j=g,f$ and $s\in\ClTm$
	\begin{IEEEeqnarray}{L}
		\PA\vdash \ClTm(\ud js) \land \val{(\ud js)} = j(\val s). \label{eq_PAf}
	\end{IEEEeqnarray}
Let $\natural$ be a translation of \ltd into \ltt that preserves the arithmetical vocabulary, commutes with connectives and quantifiers, and is so that
	\begin{IEEEeqnarray*}{RCL+C+RCL}
		\T^\natural s &:=& \Tu fs &\text{ and } & \D^\natural s &:=& \Tz gs \lor\Fz gs.
	\end{IEEEeqnarray*}
\end{dfn}

As already observed for $k$ after \Cref{dfn_flat}, we have $f\ng x=\neg f x$ for all $x$, and therefore:
	\begin{IEEEeqnarray}{L}
		\F^\natural x \Lra \Fu fx. \label{eq_Fnataux}
	\end{IEEEeqnarray}
A \textit{persistence axiom} (cf. \cite[p.~220]{can96Logical}) states that whatever is declared true at a lower level is declared true at a higher level. In the case of \kft, this amounts to $\forall\vphi\in\Sent_{\lkf}\left(\Tz\vphi\ra \Tu\vphi\right)$. 
\begin{obs}\label{obs_kf2persistence}
    By formal induction on \vphi, we have $\kftcs\vdash\forall\vphi\in\Sent_{\lkf}\left( \Tz\vphi \ra \Tu\vphi \right)$.
\end{obs}

The \textit{strong} persistence axiom mentioned above also requires, intuitively, \textit{internal persistence}, in the sense that if a sentence $\vphi\in\lkf$ is declared true at level zero, then the sentence $\vphi^+$ is declared true at level one, 
where $\vphi^+$ is obtained, informally, from \vphi by recursively replacing arbitrarily many occurrences of \Tz with \Tu, starting from the outermost occurrences and proceeding inwards. Note that, for any $x\in\Fml_\ltd$, $fx$ arises from $gx$ in this way, namely by replacing some (coded) occurrences of \Tz with \Tu, from the outermost proceeding inwards. We thus define:
    \begin{IEEEeqnarray}{L}
        \forall\vphi\in\Sent_{\ltd}\left( \Tz g\vphi \ra \Tu f\vphi  \right). \ztag{\SP}
    \end{IEEEeqnarray}

\begin{lemma}\label{le_UBkf1-5}
	The $\natural$-translations of \msf{KF1}-\msf{KF7} are derivable in $\msf{KF_2}$.
\end{lemma}

\begin{proof}
	The argument is straightforward, since \T is mapped to \Tu and there is no significant interaction between $g$ and $f$. We verify \msf{KF2} 
as an example:	
	\begin{IEEEeqnarray*}{C+L}
		\T^\natural\val t \Lra \Tu f\val t
		\labeq{\eqref{eq_PAf}}\Lra
		\Tu\val{(\ud f t)}
		\labeq{\msf{KF_22}}\Lra
		\Tu\uTu\ud f t
		\labeq{\eqref{eq_gf}}\Lra
		\Tu f\ud\T t
		\Lra
		\T^\natural\ud\T t;
		\\
		[2mm]
		\F^\natural\val t 
		\labeq{\eqref{eq_Fnataux}}\Lra \Fu f\val t
		\labeq{\eqref{eq_PAf}}\Lra
		\Fu\val{(\ud f t)}
		\labeq{\msf{KF_22}}\Lra
		\Fu\uTu\ud f t
		\labeq{\eqref{eq_gf}}\Lra
		\Fu f\ud\T t
		\labeq{\eqref{eq_Fnataux}}\Lra
		\F^\natural\ud\T t. &
	\end{IEEEeqnarray*}
The argument for \msf{KF7} is the same, \emph{mutatis mutandis}, as the one provided for \Cref{le_kf_21-5inkfd}.
\end{proof}

\begin{lemma}\label{le_UBd1-4}
	The $\natural$-translations of \msf{D1}-\msf{D4} and \msf{D9}-\msf{D10} are derivable in $\kft$.
\end{lemma}

\begin{proof}
\msf{D1} is obvious, so we begin by verifying \msf{D2}. Fix $s\in\ClTm$; then:
	\begin{IEEEeqnarray*}{R}
		\D^\natural \val s 
		\labeq{\eqref{eq_PAf}} \Lra
		\Tz\val{(\ud gs)} \lor\Fz\val{(\ud gs)}
		\labeq{\msf{KF_22}}\Lra
		\Tz\uTz\ud gs \lor\Fz\uTz\ud gs
		\Lra
		\Tz g\ud\T s \lor\Fz g\ud\T s 
		\Lra
		\D^\natural\ud\T s.
	\end{IEEEeqnarray*}

\smallskip
We verify \msf{D3}. Fix $s\in\ClTm$ and observe that
	\begin{IEEEeqnarray*}{L}
		\D^\natural \val s 
		\labeq{\eqref{eq_PAf}} \Lra
		\Tz \val{(\ud gs)} \lor \Fz\val{(\ud gs)};
		\\
		[1mm]
		\D^\natural\ud\D s \Lra
		\Tz g\ud\D s \lor\Fz g\ud\D s 
		\Lra
		\Tz(\dor{\uTz\ud{g}s}{\uFz\ud{g}s}) \lor \Fz(\dor{\uTz\ud{g}s}{\uFz\ud{g}s}).
	\end{IEEEeqnarray*}
By \msf{KF_22, KF_24}, both $\Tz \val{(\ud gs)}$ and $\Fz\val{(\ud gs)}$ entail $\Tz(\dor{\uTz\ud{g}s}{\uFz\ud{g}s})$. As for the converse direction, $\Tz(\dor{\uTz\ud{g}s}{\uFz\ud{g}s})$ entails $\Tz \val{(\ud gs)} \lor\Fz\val{(\ud gs)}$ by \msf{KF_22, KF_24}, and $\Fz(\dor{\uTz\ud{g}s}{\uFz\ud{g}s})$ entails $\Fz\val{(\ud gs)}\land\Tz\val{(\ud gs)}$ by \msf{KF_22}-\msf{KF_24}.

For \msf{D4}, use \msf{KF3} and the fact that $\F x:\lra\T\ng x$. $\msf{D9}^\natural$ is derivable from $\Tz$-regularity. For \msf{D10}, the argument is similar to the one provided in \Cref{le_kf_21-5inkfd}, relying on the fact that, provably in \PA,  $\Sent_\lkft(gx)$ only if $\Sent_\lkfd(x)$. 
\end{proof}

\begin{lemma}\label{le_UBd8}
	The $\natural$-translation of \msf{D8} is provable in $\kft$.
\end{lemma}

\begin{proof}
	Fix $t\in\ClTm$. Then 
	\begin{IEEEeqnarray*}{LL}
		\D^\natural\val t 
		\labeq{\eqref{eq_PAf}}\Lra 
		\Tz\val{(\ud gt)} \lor\Fz\val{(\ud gt)}
		\labeq{\msf{KF_20}}\Lra 
        \Tu\uTz\ud{g}t \lor\Tu\uFz\ud{g}t
        \labeq{\msf{KF_23\text{-}4}}\Lra 
		\Tu(\uTz\ud{g} t \lor\uFz\ud{g}t) 
		\Lra
		\Tu f\ud\D t 
		\Lra
		\T^\natural\ud\D t.
		\\[1mm]
		\neg\D^\natural\val t 
		\labeq{\eqref{eq_PAf}}\Lra 
		\neg\Tz\val{(\ud gt)} \land\neg\Fz\val{(\ud gt)}
		\labeq{\msf{KF_20}}\Lra 
		\Fu\uTz\ud{g} t \land\Fu\uFz\ud{g}t
		\labeq{\msf{KF_23\text{-}4}}\Lra
		\Fu(\uTz\ud{g} t \lor \uFz\ud{g}t) 
		\Lra
		\Fu f\ud\D t 
		\Lra
		\F^\natural\ud\D t. & \hspace{-2mm} \qedhere
	\end{IEEEeqnarray*}
\end{proof}

\begin{lemma}\label{le_UBcons}
	$\kftcs\vdash\forall\vphi\in\Sent_{\ltd}\left(\Tu f\vphi\ra\neg\Fu f\vphi\right) = \msf{Cons}^\sharp$. 
\end{lemma}

We now turn to the axioms \msf{D5}-\msf{D7}, for which the use of \Const and \SP is essential. 

\begin{lemma}\label{le_euivalenceT0-T1}
	$\kftcssp\vdash\forall x\in\Sent_{\ltd}\left(\Tz gx\lor\Fz gx\ra \left(\Tu fx\lra\Tz gx\right)\right)$.
\end{lemma}

\begin{proof}
    Reasoning informally within \kftcssp, assume $\Tz gx\lor\Fz gx$. By \SP, it suffices to show $\Tu fx\ra\Tz gx$, so assume $\Tu fx$. By consistency, $\neg\Fu fx$. Since $\Fz gx\ra\Fu fx$ is \SP applied to $\ng x$, by contraposition we obtain $\neg\Fz gx$. Therefore, we conclude $\Tz gx$ by $\Tz gx\lor\Fz gx$.
\end{proof}

Observe that the above lemma also yields $\Tz gx\lor\Fz gx\ra \left(\Fu fx\lra\Fz gx\right)$.

\begin{lemma}\label{le_UBd5-7}
	The $\natural$-translations of \msf{D5}-\msf{D7} are derivable in \kftcssp.
\end{lemma}

\begin{proof}
We begin with \msf{D5}. First, observe the following equivalences
	\begin{IEEEeqnarray*}{L}
		\D^\natural(\and xy) \Lra \left(\Tz(\and{gx}{gy}) \lor \Fz(\and{gx}{gy})\right)
		\labeq{\msf{KF_24}}\Lra
		(\Tz gx\land \Tz gy) \lor (\Fz gx\lor \Fz gy);
		\\[1.5mm]
		\D^\natural x\land\D^\natural y \Lra (\Tz gx\lor\Fz gx) \land (\Tz gy\lor\Fz gy);
		\\
		\D^\natural x\land\F^\natural x \Lra (\Tz gx\lor\Fz gx) \land \Fu fx \labeq{\text{\Cref{le_euivalenceT0-T1}}}\Lra\Fz gx;
		\\
		\D^\natural y\land\F^\natural y \Lra (\Tz gy\lor\Fz gy) \land \Fu fy\labeq{\text{\Cref{le_euivalenceT0-T1}}}\Lra\Fz gy.
	\end{IEEEeqnarray*}
Both implications then follow by logic.

The argument for \msf{D6} is essentially the same. As for \msf{D7}, fix $x\in\Sent_\ltd$. Then 
	\begin{IEEEeqnarray*}{L+L}
		\D^\natural x\Lra\Tz gx\lor\Fz gx 
		\labeq{\SP}\Ra
		\Tu fx\lor\Fu fx 
		\labeq{\Const}\Lra
		\neg\Tu fx\lra\Fu fx 
		\Lra
		\neg\T^\natural x\lra\F^\natural x. & \qedhere
	\end{IEEEeqnarray*}
\end{proof}

By Lemmata \ref{le_UBkf1-5}, \ref{le_UBcons}, \ref{le_UBd1-4}, \ref{le_UBd5-7}, and \ref{le_UBd8}, we finally obtain:

\begin{theorem}\label{th_UBsp}
	The translation $\natural$ is a relative interpretation of $\kfdcs$ in $\kftcssp$.
\end{theorem}

\section{Reduction of \kftcssp to \kftcs}\label{sec:redKFTCSSP_to_KFTCS}

In this section, it is shown that \SP does not add arithmetical content to \kftcs. To this end, we provide an interpretation of \kftcssp in \kftcs. Let us begin with the informal idea.
The axioms \SP imposes that \Tu must say about $f\vphi$ whatever \Tz says about $g\vphi$, otherwise it's just like the \Tu of \kftcs. In \kftcs, this can be simulated by translating $\Tu y$ by cases, namely: $\Tu y$ means that either $y$ is the $f$-image of an \lkfd-sentence $x$ whose $g$-image is true$_0$, or else the $g$-image of $x$ is neither true$_0$ nor false$_0$ and (the translation of) $y$ is true$_1$. In order to express that ``$y$ is the $f$-image of an \lkfd-sentence $x$ whose $g$-image is true$_0$'', we will use an auxiliary function $\gt = g\circ f^{-1}$. That is, $\gt y$ is $gr$ if $y=fr$ for some $r\in\Fml_{\lkfd}$, and $\gt y = 0$ otherwise. Then, $\Tz\gt y$ expresses the intended meaning.

In what follows, we let $\detz(x):=\Tz x\lor\Fz x$ and $\ud\detz(x):=\dor{\uTz x}{\uFz x}$. Moreover, let us observe that $f$ is injective on $\Fml_{\lkfd}$ and its inverse $f^{-1}$ is p.r., since the clauses of \eqref{eq_gf} have pairwise disjoint images and can be read backwards. The only case that needs checking is the clause for $\ud\D y$: its output $\dor{\uTz\ud gy}{\uFz\ud gy}$ is officially the negation $\ng{(\and{\ng{\uTz\ud gy}}{\ng{\uFz\ud gy}})}$. But this is not in the image of the clause for $\ng z$, as that would require a \Tz-atom in the range of $f$.

\begin{dfn}\label{dfn:star}
Define the function $\tilde{g}$ as follows:
    \[
        \tilde{g}x :=
        \begin{cases}
            gr & \text{if } x=fr \text{ for some } r\in\Fml_{\lkfd}
            \\
            0 & \text{otherwise}
        \end{cases}
    \]
\gt is well-defined since $f$ is injective on $\Fml_{\lkfd}$, and it is p.r.~since  $f^{-1}$ is.
Then define the function $e$ by the recursion theorem as follows
    \begin{IEEEeqnarray*}{L}
        ex:=
        \begin{cases}
			x & \text{if $x\in\mrm{At}_{\lpa}$}
			\\
			\uTz y & \text{if $x = \uTz y$, for $y\in\mrm{Term}$}
			\\
			\dor{\uTz\ud{\gt}y}{\left(\and{\ng{\ud\detz}(\ud{\gt}y)}{\uTu\ud ey}\right)} & \text{if $x = \uTu y$, for $y\in\mrm{Term}$}
			\\
			\ng e y & \text{if $x=\ng y\in\mrm{Fml}_{\ltt}$}
			\\
			\and{ey}{ez} & \text{if $x=\and yz\in\mrm{Fml}_{\ltt}$}
			\\
			\all{v}{ey} &  \text{if $x=\all vy\in\mrm{Fml}_{\ltt}$}
			\\
			0 & \text{otherwise}
		\end{cases}        
    \end{IEEEeqnarray*}
For $i\in\{\gt, e\}$ we require the same clauses as those for the other translation functions, namely, for $s\in\ClTm$
    \begin{IEEEeqnarray}{L}
		\PA\vdash \ClTm(\ud is) \land \val{(\ud is)} = i(\val s). \label{eq_PAe}
	\end{IEEEeqnarray}
Let $\star$ be a translation of \ltt into itself that is the identity on $\lpa\cup\{\Tz\}$, commutes with connectives and quantifiers, and is so that
	\begin{IEEEeqnarray*}{L}
		\Tu^\star s := \Tu es.
	\end{IEEEeqnarray*}
By our convention that $\ng x=0$ for $x\notin\Fml_\lkft$, we have $e\ng x=\ng ex$ for all $x$, hence $\Fu^\star x\Lra\Fu ex$.
\end{dfn}

\begin{obs}\label{obs:ePA}
    Provably in (a fragment of) \PA:
    \begin{enumeratei}
        \item $e$ is the identity on $\Fml_{\lkf}$,
        \item $\gt(fr)=gr$ for $r\in\Fml_{\lkfd}$,
        \item $ex\in\Sent_{\lkft}$ iff $x\in\Sent_{\lkft}$ (this is because $e$ maps formulae to formulae with the same free variables and everything else to $0$).
    \end{enumeratei}
\end{obs}

\begin{lemma}\label{le:evaluation}
    Provably in \kft, for all $t\in\ClTm$:
    \begin{IEEEeqnarray*}{L}
        \Tu e\uTu t \lra \Tz\gt\val t \lor (\neg\detz(\gt\val t) \land \Tu e\val t)
        \\
        \Fu e\uTu t \lra \neg\Tz\gt\val t \land (\detz(\gt\val t) \lor \Fu e\val t).
    \end{IEEEeqnarray*}
\end{lemma}

\begin{proof}
    We show the second claim; the first is shown analogously.
    \begin{IEEEeqnarray*}{RR}
        \Fu e\uTu t = \Fu\left( \dor{\uTz\ud{\gt}t}{\left(\and{\ng{\ud\detz}{\ugt}t}{\uTu\ud et}\right)} \right) 
        \labeq{\msf{KF_23\text{-}4}}\Lra
        \Fu\uTz\ud{\gt}t \land \Fu \left(\and{\ng{\ud\detz}{\ugt}t}{\uTu\ud et}\right)
        \\
         \labeq{\msf{KF_20, KF_22\text{-}4},\eqref{eq_PAe}}\Lra
         \neg\Tz\gt\val t \land \left(\detz(\gt\val t) \lor \Fu e\val t \right) & \hfill \qedhere
    \end{IEEEeqnarray*}
\end{proof}

\begin{lemma}\label{le:agreement}
    $\kftcs\vdash\forall\vphi\in\Sent_{\lkfd}\left( (\Tz g\vphi \ra \Tu ef\vphi) \land (\Fz g\vphi \ra \Fu ef\vphi)\right)$.
\end{lemma}

\begin{proof}
    By induction on the positive complexity of \vphi. Arithmetical literals are immediate and connective and quantifiers readily follow by IH. So we focus on truth and determinateness literals. Let first $\vphi=\ud\T s$. Then 
    \begin{IEEEeqnarray*}{R}
        \Tz g\ud\T s = \Tz\uTz\ud gs \labeq{\msf{KF_22}}\Lra \Tz g\val s \labeq{\text{\Cref{obs:ePA}(ii)}}= \Tz \gt(f\val s) \Ra \Tz\gt(f\val s)\lor \bigl(\neg\detz(\gt(f \val s))\land\Tu e(f\val s)\bigr)
        \\
        \labeq{\text{\Cref{le:evaluation}},\eqref{eq_PAe}}\Lra
        \Tu e\uTu\ud fs = \Tu ef\ud\T s
    \end{IEEEeqnarray*}
For $\Fz g\ud\T s \ra \Fu ef\ud\T s$ we use \Const:
    \begin{IEEEeqnarray*}{R}
        \Fz g\ud\T s = \Fz\uTz\ud gs \labeq{\msf{KF_22}}\Lra \Fz g\val s 
        \labeq{\text{\Cref{obs:ePA}(ii)}, \msf{KF_27}}\Lra \Fz \gt(f\val s) 
        \labeq{\Const}\Lra \neg\Tz \gt(f\val s) \land \detz(\gt(f\val s))
        \\
        \Ra \neg\Tz \gt(f\val s) \land \left(\detz(\gt(f\val s)) \lor \Fu e(f\val s)\right)
        \labeq{\text{\Cref{le:evaluation}}, \eqref{eq_PAe}}\Lra \Fu e\uTu\ud fs = \Fu ef\ud\Ts.
    \end{IEEEeqnarray*}

The argument is dual for $\ng{\ud\T}s$.

Let $\vphi=(\ng)\ud\D s$. Then $f\vphi=g\vphi\in\lkf$, hence $e(f\vphi)=g\vphi$ by \Cref{obs:ePA}(i), so we are done by \Cref{obs_kf2persistence}. 
\end{proof}

\begin{lemma}\label{le_kftcssp_in_kftcs}
    The translation $\star$ is a relative interpretation of \kftcssp in \kftcs.
\end{lemma}

\begin{proof}
    The axioms for \Tz are obvious since $\star$ only changes \Tu. The axioms for \Tu are mostly straightforward, except for \msf{KF_22}, for which we reason as follows. First, observe that
    \begin{IEEEeqnarray}{L}
        \Tu^\star\uTu t = \Tu e\uTu t \labeq{\text{\Cref{le:evaluation}}}\Lra \Tz\gt\val t \lor (\neg\detz(\gt\val t) \land \Tu e\val t) \label{eq_bullet}
        \\
        \Tu^\star\val t = \Tu e\val t \nonumber
    \end{IEEEeqnarray}
If $\neg\detz(\gt\val t)$, then in particular $\neg \Tz\gt\val t$, and thus $\Tu^\star\uTu t\lra\Tu e\val t\lra\Tu^\star\val t$. If $\detz(\gt\val t)$, then by \msf{KF_27} $\gt\val t$ is a sentence, and thus $\gt\val t\neq0$. The definition of \gt gives $r\in\Fml_{\lkfd}$ such that $\val t=fr$ and $\gt\val t = gr$. Since $\gt\val t=gr$, we obtain $r\in\Sent_{\lkfd}$. Then we have 
    \begin{IEEEeqnarray*}{L}
        \Tu e(\uTu t) \labeq{\eqref{eq_bullet}}\Lra\Tz\gt\val t = \Tz gr 
        \labeq{\Const, \text{\Cref{le:agreement}}}\Lra\Tu efr = \Tu e\val t.
    \end{IEEEeqnarray*}
The argument showing $(\Fu\val t\lra\Fu\uTu t)^\star$ is essentially the same.

$\msf{Cons_2^\star}$ follows from \Const. Finally, note that $\SP^\star = \forall\vphi\in\Sent_\lkfd\left( \Tz g\vphi \ra \Tu ef\vphi \right)$, so we are done by \Cref{le:agreement}.
\end{proof}

Since $\star$ is identity on $\lpa\cup\{\Tz\}$, we obtain

\begin{prop}\label{pr_kftcssp_in_kftcs}
    \kftcssp is \lkf-conservative over \kftcs.
\end{prop}

Taken together, \Cref{pr_kftcssp_in_kftcs} and Theorems \ref{th_LB}, \ref{th_FujCon}, and \ref{th_UBsp} yield

\begin{theorem}
	$\vert\kfd\vert = \vert\kfdcs\vert = \vert\kft\vert = \vert\widehat{\msf{ID}}_\msf{2}\vert = \alpha_2$.
\end{theorem}

\section{Theories with restricted induction}\label{sec:add_res}
By inspecting the proof of the lower bound for \kfd, it can be observed that the induction schema is used only to interpret the induction schema of \kft. It is then straightforward to verify that restricting induction on both sides preserves the interpretations from §§\ref{sec:LB}-\ref{sec:UB}. Let us provide some details (using a sloppy notation).

\begin{dfn}
    For $i=0,1$, let $\msf{tot}(\vphi, \T_i):=\forall x(\T_i\vphi(x)\lor\F_i\vphi(x))$ and define
    \begin{IEEEeqnarray}{L+L}
		\forall\vphi\in\Fml^1_\Lmc\left(\big(\T_i\vphi(0) \land \forall x \left( \T_i\vphi(x) \ra \T_i\vphi(x+1) \right)\big) \ra \forall y \T_i\vphi(y)\right) & \ztag{$\msf{IND}^\msf{int}_{\Lmc,\T_i}$}
        \\
        \forall\vphi\in\Fml^1_\Lmc \left(\msf{tot}(\vphi,\T_i) \ra \big(\T_i\vphi(0) \land \forall x \left( \T_i\vphi(x) \ra \T_i\vphi(x+1) \right)\big) \ra \forall y \T_i\vphi(y)\right) & \ztag{$\msf{IND}^\msf{tot}_{\Lmc,\T_i}$}
	\end{IEEEeqnarray}
\end{dfn}

Since \Lmc will always be clear from context, we just write $\msf{IND}^\msf{int}_{\T_i}$ and $\msf{IND}^\msf{tot}_{\T_i}$.

\begin{dfn}
\hfill
    \begin{itemize}
        \item[-] $\kftint \, (+\Const)$ and $\kfttot \ (+\Const)$ are the systems obtained from $\kft \ (+\Const)$ by replacing $\msf{IND}_{\ltt}$ with $\msf{IND}^\msf{int}_{\Tu}$ and, respectively, $\msf{IND}^\msf{tot}_{\Tu}$.
         \item[-] $\kfdint \ (\kfdintcs)$ and $\kfdtot \, (\kfdtotcs)$ are the systems obtained from $\kfd \, (\kfdcs)$ by replacing $\msf{IND}_{\lkfd}$ with $\msf{IND}^\msf{int}_{\T}$ and, respectively, $\msf{IND}^\msf{tot}_{\T}$.\footnote{
            Recall that, for \kf and \kfd, we write \T instead of \Tz.
        }
    \end{itemize}
\end{dfn}

\begin{lemma}
    The translation $\flat$ from \Cref{dfn_flat} is a relative interpretation of \kftr (\kfttot, \kftint) in \kfdr (\kfdtot, \kfdint).
\end{lemma}

\begin{proof}
    In light of \Cref{th_LB}, it suffices to observe that the $\flat$-images of $\msf{IND}^\msf{int}_{\Tu}$, $\msf{IND}^\msf{tot}_{\Tu}$ follow from, respectively, $\msf{IND}^\msf{int}_{\T}$ and $\msf{IND}^\msf{tot}_{\T}$.
\end{proof}

Similarly, by noticing that the $\natural$-translations of induction schemata are derivable in the corresponding \kft-variant, the arguments from §\ref{sec:UB} readily entail that

\begin{lemma}
	The translation $\natural$ from \Cref{dfn_natural} is a relative interpretation of \kfdrcs (\kfdtotcs, \kfdintcs)  in \kftrcssp (\kfttotcssp, \kftintcssp).
\end{lemma}

We now provide bounds for induction-variants of \kfdcs via an analysis of the corresponding variants of \kft, showing in particular that \Const and \SP can be conservatively added to any induction variant of \kft.

\subsection{Arithmetical and total induction}

The conservativity of \kfdrcs over \PA has already been established in §\ref{sec:semantics}, so we consider \kfttotcs.\footnote{
    Let us observe that \kftrcssp can be shown to be conservative over \PA by a simple model-theoretic argument. Given a model $\Nmc\vDash\PA$, let $K_0$ be the Cantini-Kripke Jump for \Tz over \Nmc, and $K_1$ be the Cantini-Kripke Jump for \Tu over $(\Nmc, X_0)$, where $X_0$ is an interpretation for \Tz. Then the structure $(\Nmc, S_0, S_1)$ where $S_0$ is the minimal fixed-point of $K_0$ over \Nmc and $S_1$ is the minimal fixed-point of $K_1$ over $(\Nmc, S_0)$ is a model of \kftrcssp. \SP can be verified reasoning by induction on the $K_0$-rank of $ga$.
}

\begin{lemma}
    \kfttotcssp is \lpa-conservative over \kfcs.
\end{lemma}

\begin{proof}[Proof Sketch]
    Suppose $\kfcs\not\vdash\sigma$ for $\sigma\in\lpa$, and let $(\Nmc, S_0)\vDash\kfcs$ be a countable recursively saturated model such that $(\Nmc, S_0)\vDash\neg\sigma$. Let $K_1$ be the Cantini-Kripke Jump, or \CK-Jump for short, for \Tu defined over \Nmc with $S_0$ as parameter for \Tz.\footnote{
        The \CK-Jump can be defined in the usual way. See e.g. the formula $\Acr_1(x, X, Y)$ below.
    }
Then, for $B:=(S_0\cap\Sent_\lkf) \cup \text{ the closure under regularity of }\{fa\mid a\in\Sent_{\lkfd} \ \& \ ga\in S_0\}$,\footnote{
    Write $x\approx y$ if $x$ and $y$ coincide once every closed term in them is replaced by the numeral of its value; this relation is primitive recursive. The closure of a set $X$ under regularity is $X^*:=\{x\mid x \approx y\text{ for some }y\in X\}$. 
} 
define $S^0_1:=B$, $S^{n+1}_1:=K_1(S^n_1)$ and $S^\omega_1:=\bigcup_{n<\omega}S^n_1$. We claim:
    \begin{enumeratei}
        \item $S_1^k\subseteq S_1^j$ for every $k\leq j$ and every $S_1^i$ is consistent and regular;
        \item $S^\omega_1$ is a fixed-point of $K_1$;
        \item $(\Nmc, S_0, S^\omega_1)\vDash\kfttotcssp$.
    \end{enumeratei}
By monotonicity of $K_1$, for claim (i) it suffices to show that $B\subseteq K_1(B)$, reasoning by cases in a straightforward manner. For example, using a sloppy notation, if $\vphi\equiv\neg\Tz s\in S_0$, then $\neg\val s \in S_0$; by consistency of $S_0$, $\val s \notin S_0$ hence $\neg\T_0 s\in K_1(B)$. If $\vphi\equiv\psi\land\xi\in S_0$, then $\psi,\xi\in S_0$ hence $\vphi\in K_1(B)$. Similarly, if $\vphi=fa$, then one reasons by cases depending on the shape of $a\in\Sent_\lkfd$. If $a=\D s$, then $fa=ga\in S_0\subseteq K_1(B)$. If $a=\T s$, then $ga=\Tz gs\in S_0$, which is the case iff $gs\in S_0$, and hence $fs\in B$, which in turn entails $\Tu fs=f(\Ts)\in K_1(B)$. 
The consistency and regularity of each $S_1^i$ then follows from the fact that $K_1$ preserves regularity (see \cite[Proposition 5.4(i)]{can89Notes}) and that $K_1(X)$ is consistent whenever $X\subseteq K_1(X)$ is consistent; and it can be shown that $B$ is both consistent and regular.

Claims (ii) and (iii) are essentially due to \cite[pp.~198-200]{clw17models}. Since they require a slight modification in order to be used in the present context, we quickly summarize the main ingredients of the argument for the sake of the reader's convenience. For (ii), by monotonicity of $K_1$ it suffices to show that $K_1(S^\omega_1)\subseteq S^\omega_1$. Every clause of the jump is finitary except for the clause concerning $\forall x\vphi$. Thus, the key fact to show is that, if $\vphi(t)\in S^\omega_1$ for all closed terms $t$ (in the sense of \Nmc), then $\forall x\vphi(x)\in S^\omega_1$. So assume $\vphi(t)\in S^\omega_1$ for all such $t$. Towards a contradiction, assume there is no single $n$ such that $\vphi(t)\in S^n_1$ for all $t$. Let $J^n\in\lkf$ be a formula defining $S^n_1$ so that $(\Nmc, S_0)\vDash J^n(\vphi)$ iff $\vphi\in S^n_1$,\footnote{
    We can take $J^0(x) := \exists y\left(x\approx y\land\left((\Tz y\land\Sent_{\lkf}(y))\lor\exists a\in\Sent_{\lkfd}\left(y=fa \land\Tz ga\right)\right)\right)$, and $J^{n+1}(x)$ to be defined in the standard way.
} 
and consider the following recursive set of formulae
    \[
        p(x):=\{ x\in\ClTm\land\neg J^n(\vphi(x)) \mid n\in\omega \}.
    \]
Since $(\Nmc, S_0)\vDash\forall x\left( x\in J^k\ra x\in J^i\right)$ for all $k\leq i$, and since we are assuming that there is no single $n$ such that $\vphi(t)\in S^n_1$ for all $t$, it follows that for every $n$ there is a $t$ such that $(\Nmc, S_0)\vDash\neg J^n(\vphi(t))$. Hence, $p(x)$ is finitely realizable in $(\Nmc, S_0)$, and thus by recursive saturation there is some $s$ realizing it. But this entails that $\vphi(s)\notin S^\omega_1$, contradicting our assumption. Hence, there is an $n$ such that $\vphi(t)\in S^n_1$ for all $t$, and hence $\forall x\vphi x\in S_1^{n+1}$.

For (iii), it is clear that $(\Nmc, S_0, S^\omega_1)\vDash\kftrcssp$, since $S_0$, $S^\omega_1$ are both fixed-points and $B\subseteq S^\omega_1$; so we only verify total induction. Let $\vphi(x)$ be an arbitrary formula (in the sense of \Nmc) such that $(\Nmc, S_0, S^\omega_1)\vDash\forall t\left( \Tu\vphi(t) \lor \Fu\vphi(t) \right)$. By compositionality, $(\Nmc, S_0, S^\omega_1)\vDash\Tu\left(\forall x\left(\vphi(x) \lor \neg\vphi(x) \right)\right)$. Thus, there exists an $n$ such that, for all closed terms $t$, $\vphi(t)\in S^n_1$ or $\neg\vphi(t) \in S^n_1$. 
Since $S^n_1\subseteq S^\omega_1$, by consistency of $S^\omega_1$ we obtain $(\Nmc, S_0, S^\omega_1)\vDash\Tu\vphi(t)\lra J^n(\vphi(t))$. In fact, $J^n(\vphi(t))\ra\Tu\vphi(t)$ follows from $S^n_1\subseteq S^\omega_1$; for the converse, if $\vphi(t)\in S^\omega_1$ but $\vphi(t)\notin S^n_1$, then we would have $\neg\vphi(t)\in S^n_1\subseteq S^\omega_1$. Finally, since $J^n(\vphi(x))\in\lkf$, we have full induction for it in $(\Nmc, S_0)$, thus we obtain induction for $\Tu\vphi(x)$ in $(\Nmc, S_0, S^\omega_1)$.
\end{proof}

Since full disquotation of \Tu over \Tz makes every \lkf formula \Tu-total in \kfttot, we have: 

\begin{lemma}
    \kf is a subtheory of \kfttot.
\end{lemma}

Thus, since $\vert\kf\vert = \vert\kfcs\vert$, we obtain:

\begin{prop}
    $\vert\kfdtot\vert = \vert\kfdtotcs\vert = \vert\kfttot\vert = \vert\kfttotcssp\vert = \vert\kfcs\vert = \vert\kf\vert = \vphi_{\varepsilon_0}0$.
\end{prop}

\subsection{Internal induction}
The analysis of theories with internal induction requires more work, but the arguments mostly rely on known results, so we keep the presentation sketchy.

\subsubsection{Lower bound}
Beginning with the lower bound, we prove $\vphi_{\vphi_\omega0}0\leq\vert\kftint\vert$. By \Cref{le:pos+negIND} below, the \Tu-level of \kftint is, over \lkf, a theory with positive induction in the sense of \cite{js96Some}. J\"ager and Strahm show that the fixed-point theory $\msf{ID}^\#_1$, in which induction is restricted to formulae positive in the fixed-point predicates, has ordinal $\vphi_\omega0$ \cite[Theorem~38]{js96Some}.\footnote{
    They also note that the same holds for $\msf{Ref}^\#(\PA)$, i.e., Feferman's $\msf{Ref}(\PA)$ \cite{fef91Reflecting} with induction restricted to formulae positive in \T and \F \cite[§5]{js96Some}.
}
The present treatment of the \Tu-level adapts their lower-bound argument \cite[§2.2]{js96Some}, except that what we embed in \kftint is not $(\Pi^0_1\text{-}\msf{CA})_{<\omega^\omega}$, but rather the system $\RT{<\omega^\omega}[\lkf]$ of ramified truth over \lkf up to any ordinal below $\omega^\omega$.
We use a standard notation system for ordinals (see e.g.~\cite{poh09Proof}), and we define the ramified systems of truth over two different base languages, namely \lpa and \lkf. We use $\Lmc_B$ as a variable ranging over \lpa and \lkf.

\begin{dfn}
For $\alpha<\Gamma_0$, let $\Lmc^B_\alpha$ be the language $\Lmc_B\cup\{\T_\beta\}_{\beta<\alpha}$, for $\T_\beta$ a new unary predicate.\footnote{
    We use the predicates \Tz and \Tu for both the ramified language and \lkft. Context will always make clear which one we are referring to.
}
Then $\RT{<\alpha}[{\Lmc_B}]$ is the theory over $\Lmc^B_\alpha$ expanding $\PA_{\Lmc^B_\alpha}$ with the following axioms, for each $\gamma<\beta<\alpha$, and where $\Sent_{\beta}, \Fml^1_{\beta}$ stand for $\Sent_{\Lmc^B_\beta}, \Fml^1_{\Lmc^B_\beta}$:
\begin{enumerate}[label=$\RT{}\arabic*_{<\alpha} $] 
    \item $\forall s\, \T_\beta(\ud Rs) \lra R\val s$, for each $\Lmc_B$-atomic $R$.
    \item $\forall \vphi \in \Sent_{\beta}$ $\T_\beta(\neg \vphi) \lra \neg \T_\beta(\vphi).$
    \item $\forall \vphi, \psi \in \Sent_{\beta}$ $\T_\beta(\vphi \land \psi) \lra (\T_\beta(\vphi) \land \T_\beta(\psi)).$
    \item $\forall \vphi(v) \in \Fml^1_{\beta}$ $\T_\beta(\forall v \vphi(v)) \lra \forall s\, \T_\beta(\vphi(s)).$
    \item $\forall s \left(\val s \in \Sent_{\gamma} \ra \bigl(\T_\beta\ud\T_\gamma
s \lra \T_\gamma(\val s)\bigr)\right)$
    \item $\forall s\forall\delta\prec\beta \left(\val s \in \Sent_{\delta} \ra \bigl(\T_\beta \ud\T_\delta s \lra \T_\beta(\val s)\bigr)\right)$.
\end{enumerate}
\end{dfn}

The proof of the lower bound consists of the following main steps. We first observe that \kftint truth-defines $\RT{<\omega^\omega}[\lkf]$. This follows from the full transparency of \Tu over \lkf along with the fact that \kftint derives transfinite induction up to any ordinal below $\omega^\omega$ for a certain fragment of \ltt. By the truth-definability of $\RT{<\omega^\omega}[\lkf]$, we obtain that $\kftint\vdash\TI(<\vphi_\omega0, \lkf)$. Finally, by \cite[Lemma 36]{fuj10Relative} this is shown to yield the truth-definability of $\RT{<\vphi_\omega0}[\lpa]$ in \kftint, and thus the derivability of $\TI(<\vphi_{\vphi_\omega0}0, \lpa)$.

To make this formally precise, we define a translation of $\Lmc^\kf_\alpha$ in \ltt. We preliminarily define a binary p.r.\@ function $q$ such that 
\begin{align*} 
q ( \beta, x ) := 
\begin{cases} 
x & \text{if $x \in \mathrm{Sent}_{\beta}$} 
\\ 
\ulcorner 0 = 1 \urcorner & \text{otherwise}. 
\end{cases} 
\end{align*} 
We next define a unary p.r.\@ function $p$ on (the codes of) $\Lmc^\kf_{\alpha}$-sentences $z$
so that:

\[ 
p ( z ) := 
\begin{cases} 
z & \text{if $z$ is an atomic $\lkf$-sentence}
\\ 
\uTu \bigl( \ud p\ud q ( \gamma, s ) \bigr) 
& \text{if $z = \ud\T_{\gamma} s$ for $s \in \ClTm$} \\ 
\ng{\uTu}\ud p w 
& \text{if $z = \ng w$} \\ 
\and{\uTu\ud pw}{\uTu\ud pu}
& \text{if $z = \and{w}{u}$} \\ 
\all v \uTu\ud p(w\dot v) & \text{if $z = \all vw$} 
\\
\cor{0=1} & \text{otherwise}
\end{cases}
\]
Then, to obtain the truth-definition of $\RT{<\omega^\omega}[\lkf]$, the formula $\T_\beta x$ is translated as $\theta_\beta(x)$, where $\theta_\beta(x):=\Tu pq(\beta, x)$. Moreover, we define
    \begin{IEEEeqnarray*}{L}
        \Tot(\beta):=\forall x(\Sent_\beta(x)\to\T_1 px\lor\F_1 px),
        \\
        \Con(\beta):=\forall x(\Sent_\beta(x)\to\neg(\T_1 px\land\F_1 px)).
    \end{IEEEeqnarray*}
Let $\Pos$ ($\Neg$) be the class of $\ltt$-formulae in which $\T_1$ occurs only positively (negatively). Note that $\Tot\in\Pos$, while $\Con\in\Neg$. 

\begin{lemma}\label{le:commutation_kftr}
    $\kftr\vdash\forall\beta\,\forall y\,\left(\Tot(\beta)\land\Con(\beta)\to \bigl(\neg\Tu pq(\beta,y)\lra\Fu pq(\beta,y)\bigr)\right)$. In particular, if $\kftr\subseteq S\vdash\Tot(\beta)\land\Con(\beta)$, then $S\vdash\neg\theta_\beta(t)\lra\Fu pq(\beta,t)$ for every term $t$.
\end{lemma}

\begin{proof}
We reason informally in $\kftr$. Fix $\beta,y$ and assume $\Tot(\beta)\land\Con(\beta)$. We distinguish two cases. Note that the case distinction is on the $\lpa$-formula $\Sent_\beta(y)$ and involves no induction. 
\begin{itemize}
  \item \emph{$\Sent_\beta(y)$.} Then $q(\beta,y)=y$. $\Tot(\beta)$ gives
    $\Tu p(y)\lor\Fu p(y)$, and $\Con(\beta)$ gives $\neg(\Tu p(y)\land\Fu p(y))$. Together they give $\neg\Tu p(y)\lra\Fu p(y)$.
  \item \emph{$\neg\Sent_\beta(y)$.} Then $q(\beta,y)=\cor{0=1}$, and $p(\cor{0=1})=\cor{0=1}$. By $\mathsf{KF}_21$, $\neg\T_1\cor{0=1}\land\Fu\cor{0=1}$, so both sides of the equivalence are true.
\end{itemize}
For the `in particular', take $y:=t$ and recall that $\theta_\beta(t)$ is
$\Tu pq(\beta,t)$.
\end{proof}

\begin{lemma}\label{le:pos+negIND}
    \kftint proves induction for every $\vphi\in\Pos\cup\Neg$. Thus, since $\lkf\subset\Pos$, \kftint proves induction for \lkf.
\end{lemma}

\begin{proof}
    For $\vphi\in\Pos$, one shows by meta-induction on \vphi that $\kftint\vdash\Tu\cor{\vphi(\dot x)}\lra\vphi(x)$ (cf.~\cite[Prop.~3]{jkss99Proof}); for $\vphi\in\Neg$, use the fact that \kftint has induction for $\neg\Tu$ by Parsons' trick (cf. \cite{hh06Axiomatizing}). 
\end{proof}

The following is also readily established:

\begin{lemma}[Cf.~{\cite[Lemma~2, Corollary~3]{js96Some}}]\label{le:pos+negTI}
    For all $n$ and all $\vphi\in\Pos\cup\Neg$, $\kftint\vdash\TI(\omega^n, \vphi)$.
\end{lemma}

\begin{proof}[Proof Sketch]
    It suffices to show, by meta-induction on $k$, that the following jump formula is derivable for any progressive $\vphi\in\Pos\cup\Neg$: 
    \[
        \forall\xi\left(\forall\eta\prec\xi\vphi(\eta) \ra \forall\eta\prec\xi+\omega^k\vphi(\eta) \right).
    \]
See \cite[Lemma 3.3]{can89Notes} and \cite[Theorem 36]{hh06Axiomatizing} for details.
\end{proof}

\begin{lemma}\label{le:ProgTot+Con}
    $\kftint\vdash\Prog(\Tot(\beta))\land\Prog(\Con(\beta))$.
\end{lemma}

\begin{proof}[Proof Sketch]
Since $\Tot(\beta)\in\Pos$ and $\Con(\beta)\in\Neg$, by \Cref{le:pos+negIND}, we can reason by induction on the complexity of codes of formulae. 
So assume $\forall\gamma\prec\beta\,\Tot(\gamma)$ to show $\forall x\left(\Sent_\beta(x) \ra \Tu px\lor\Fu px \right)$. We only show two examples and refer the reader to \cite[p.~333]{fuj10Relative} for further details.
For $x$ an atomic sentence of $\Lmc_0=\lkf$, the claim follows immediately by \msf{KF_20}-\msf{KF_21} since $p(x)=x$. Let $x=\T_{\gamma} s$ with $\gamma \prec \beta$. Then for $\Tot(\beta)$ we have to show $\Tu pq ( \gamma, s ) \vee \Fu pq ( \gamma, s )$. 
If $s \not\in \Sent_{\gamma}$, then $p \circ q ( \gamma, s ) = \cor{0 = 1}$ and thus $\Fu p \circ q ( \gamma, s )$ by \msf{KF_21}. Otherwise, $s \in \Sent_{\gamma}$ and thus $\Tu p ( s ) \vee \Fu p ( s )$ by IH on $\gamma$. Similarly for $\Con(\beta)$---note that \Const is not needed.
\end{proof}

\begin{lemma}[Cf.~{\cite[Lemma~4]{js96Some}}]\label{le:TI_Tot_Con}
    For all $n$, $\kftint\vdash\forall\beta\prec\omega^n\left(\Tot(\beta) \land \Con(\beta) \right)$.
\end{lemma}

\begin{proof}
    By \Cref{le:pos+negTI} and \Cref{le:ProgTot+Con}, since $\Tot(\beta)\in\Pos$ and $\Con(\beta)\in\Neg$.
\end{proof}

\begin{prop}[Cf.~{\cite[Theorem~6]{js96Some}}]\label{pr:RTinKFtint}
    The map $\T_\beta(x) \mapsto \theta_\beta(x)$ is a truth-definition of $\RT{<\omega^\omega}[\lkf]$.
\end{prop}

\begin{proof}
    Axioms $\msf{RT}1_{<\omega^\omega}$ and $\msf{RT}3_{<\omega^\omega}$-$\msf{RT}6_{<\omega^\omega}$ go through as in \cite[p.~333]{fuj10Relative}. $\msf{RT}2_{<\omega^\omega}$ and instances of induction need $\neg\theta_\beta(x)\lra\Fu pq(\beta,x)$. \Cref{le:commutation_kftr} yields this if $\Tot(\beta) \land \Con(\beta)$, and  the latter follows from \Cref{le:TI_Tot_Con}. Then, for $\msf{RT}2_{<\omega^\omega}$ we have $\theta_\beta(\ng x)\lra\Fu pq(\beta,x)$ by \msf{KF_22}. For induction, let \vphi be arbitrary and let $\vphi'$ be its image under $\T_\beta(x) \mapsto \theta_\beta(x)$ in negation normal form. Then every occurrence of \Tu in $\vphi'$ is of the form $\theta_\beta(x)$, i.e., $\Tu pq(\beta,x)$. Since $\neg\theta_\beta(x)\lra\Fu pq(\beta,x)$, it follows that $\vphi'$ is equivalent to a positive formula, thus we are done by \Cref{le:pos+negIND}.
\end{proof}

\begin{theorem}
    $\kftint\vdash\TI(<\vphi_{\vphi_\omega0}0, \lpa)$.
\end{theorem}

\begin{proof}
\kftint satisfies conditions (i)-(viii) of \cite[Lemma 36]{fuj10Relative} for $\T=\Tz$ and derives induction for \lkf by \Cref{le:pos+negIND}. The well-known equivalence between systems of ramified analysis and those of ramified truth entails that $\RT{<\omega^\omega}[\lkf]\vdash\TI(<\vphi_\omega0, \lkf)$. Thus, by \Cref{pr:RTinKFtint}, $\kftint\vdash\TI(<\vphi_\omega0, \lkf)$. It follows that \kftint truth-defines $\RT{<\vphi_\omega0}[\lpa]$. Since $\RT{<\vphi_\omega0}[\lpa]\vdash\TI(<\vphi_{\vphi_\omega0}0, \lpa)$, we obtain the desired conclusion.
\end{proof}

\subsubsection{Upper bound}\label{sssec:ub}
For the upper bound, we prove that \kftintcssp is arithmetically conservative over $\PA+\TI(<\vphi_{\vphi_\omega0}0)$. First note that the arguments from §\ref{sec:redKFTCSSP_to_KFTCS} go through with internal induction only. The only
inductions there, in \Cref{obs_kf2persistence} and \Cref{le:agreement}, are on
formulae in \Pos, so \Cref{le:pos+negIND} applies. 
instances of $\msf{IND}^\msf{int}_{\Tu}$. 
Thus:

\begin{lemma}
    The translation $\star$ from \Cref{dfn:star} is an interpretation of \kftintcssp in \kftintcs.
\end{lemma}

The remainder of this section is devoted to proving the following: 

\begin{theorem}\label{th:UB-kftintcs}
    If $\kftintcs\vdash\sigma$ for $\sigma\in\lpa$, then $\PA+\TI(<\vphi_{\vphi_\omega0}0)\vdash\sigma$.
\end{theorem}

In proving the theorem, we adapt the argument of \cite[§4]{fuj11Autonomous} (see \Cref{th_FujCon} above) in order to establish the conservativity of \Const over \kftint. We assume the reader is familiar with cut-elimination and asymmetric interpretation techniques employed in the proof-theoretic analysis of the systems \idh\alpha (see, e.g., \cite{can85Note,jkss99Proof}). For positive induction, these techniques were worked out in \cite[§4.3]{js96Some}: the asymmetric interpretation only needs stages below $\omega^\omega$, and predicative cut-elimination then gives bounds below $\vphi_\omega0$ \cite[Theorems~35--36]{js96Some}.

We introduce the semi-formal system \h{2} (cf.~\cite{can85Note,jkss99Proof}), formulated in a Tait-style language $\Lcr_{\h2}$ extending \lpa with unary symbols $\p0{}, \p1{}$.
We first define two inductive operator forms $\Acr_0(x, X)$ and $\Acr_1(x, X, Y)$ and fix them throughout the subsequent arguments. 
The operator form $\Acr_0(x, X)$ corresponds to the \CK-Jump for $\T_0$, and it is given by
	\begin{IEEEeqnarray*}{LL}
		\Sent_{\lkf}(x) & \land \forall yz\in\ClTm\left( x = \eqn yz \ra \val y = \val z \right)
		\\
		&\land \forall yz\in\ClTm\left( x = \ng{(\eqn yz)} \ra \val y \neq \val z \right)
		\\
		&\land \forall y\in\ClTm\left(x = \uTz y \ra \val y \in X \right)
		\\
		&\land \forall y\in\ClTm\left(x = \ng{\uTz}y \ra \ng{\val y} \in X \right)
		\\
		&\land \forall y\in\Sent_{\lkf}\left(x = \ng{\ng y} \ra y \in X \right)
		\\
		&\land \forall y,z\in\Sent_{\lkf}\left(x = \and yz \ra y \in X \land z \in X \right)
		\\
		&\land \forall y,z\in\Sent_{\lkf}\left(x = \ng{(\and yz)} \ra \ng y \in X \lor \ng z \in X \right)
		\\
		&\land \forall y,z\left(x = \all yz\in\Sent_{\lkf} \ra \forall u(z(u) \in X) \right)
		\\
		&\land \forall y,z\left(x = \ng{\all yz}\in\Sent_{\lkf} \ra \exists u(\ng z(u) \in X) \right).
	\end{IEEEeqnarray*}
The operator form $\Acr_1(x, X, Y)$ corresponds to the \CK-Jump for $\T_1$ and it is defined thus:
    \begin{IEEEeqnarray*}{LL}
		\Sent_{\lkft}(x) & \land \forall yz\in\ClTm\left( x = \eqn yz \ra \val y = \val z \right) 
		\\
		&\land \forall yz\in\ClTm\left( x = \ng{(\eqn yz)} \ra \val y \neq \val z \right)
		\\
		&\land \forall y\in\ClTm\left(x = \uTz y \ra \val y\in X \right)
		\\
		&\land \forall y\in\ClTm\left(x = \ng{\uTz y} \ra \val y\notin X \right)
		\\
		&\land \forall y\in\ClTm\left(x = \uTu y \ra \val y \in Y \right)
		\\
		&\land \forall y\in\ClTm\left(x = \ng{\uTu}y \ra \ng{\val y} \in Y\right)
		\\
		&\land \forall y\in\Sent_{\lkft}\left(x = \ng{\ng y} \ra y \in Y \right)
		\\
		&\land \forall y,z\in\Sent_{\lkft}\left(x = \and yz \ra y \in Y \land z \in Y \right)
		\\
		&\land \forall y,z\in\Sent_{\lkft}\left(x = \ng{(\and yz)} \ra \ng y \in Y \lor \ng z \in Y \right)
		\\
		&\land \forall y,z\left(x = \all yz\in\Sent_{\lkft} \ra \forall u(z(u) \in Y) \right) 
		\\
		&\land \forall y,z\left(x = \ng{\all yz}\in\Sent_{\lkft} \ra \exists u(\ng z(u) \in Y)\right).
	\end{IEEEeqnarray*}
The rules of \h{2} are those of $\omega$-arithmetic along with the following axioms and fixed-point rules for \p{i}{}.
To formulate rules for \p0{} and \p1{} uniformly, we write $\Acr_1(s, \p1{})$ instead of $\Acr_1(s, \p0{}, \p1{})$:

\smallskip
\noindent
\textit{Axioms.}
	\begin{enumerate}[label=Ax\arabic*]
	\setcounter{enumi}{-1}
		\item $\Gamma, \vphi$ for a true \lpa-literal $\vphi$.
		\item $\Gamma, \neg\p{i}{}(t), \p{i}{}(s)$ for $s^\Nbl=t^\Nbl$ and $i=0,1$.
	\end{enumerate}

\noindent
\textit{Fixed-point rules for \p{i}{}, $i=0,1$}:
	\begin{IEEEeqnarray*}{L}
	\def\fCenter{ \mbox{$\Gamma$} }
		\ax{ \fCenter, \Acr_i(s, \p{i}{})}
        \rl{$\msf{FP}_i^+$}
			\un{ \fCenter, \p{i}{}(s) }
		\Dis
        \hspace{20mm}
        \ax{ \fCenter, \neg\Acr_i(s, \p{i}{})}
       \rl{$\msf{FP}_i^-$}
			\un{ \fCenter, \neg\p{i}{}(s) }
        \Dis
	\end{IEEEeqnarray*}
\h1 is \h2 without \p1{} and its fixed-point rules; \h0 is $\omega$-arithmetic. We now show that, if $\kftint\vdash\vphi(x)$, then $\vphi(\ovl n)[\p{}{}/\T]$ is derivable in \h2, where $\vphi(\ovl n)[\p{}{}/\T]$ is obtained from $\vphi(\ovl n)$ by replacing every occurrence of $\T_i$ with $\p{i}{}$; we let $\Gamman[\p{}{}/\T]:=\{\vphi(\ovl n)[\p{}{}/\T]\mid\vphi(x)\in\Gamma\}$. To state this precisely, we assume that \kftint is given as a Tait-calculus allowing for partial cut-elimination, i.e., every derivation of \kftint can be effectively transformed into one in which all cuts are on atomic formulae (see, e.g., \cite[§9]{can89Notes}). Moreover, we let $\msf{H}_{k}\svdash{\star}{\delta}\Gamma$ (respectively $\msf{H}_{k}\svdash{\star}{<\delta}\Gamma$) mean that $\Gamma$ is derivable in $\msf{H}_{k}$ with a derivation of height $\leq\delta$ (resp.~$<\delta$) and cuts restricted to formulae of the form $\p{i}{}, \neg\p{i}{}$.

\begin{lemma}\label{le:embedding_kftint-h2}
    If $\kftint\vdash\Gamma(\vec x)$, then there exists $k<\omega$ such that $\h{2}\svdash{\star}{<\omega\cdot k}\Gamman[\p{}{}/\T]$, for all $\vec n$ with $k$ independent of $\vec n$.
\end{lemma}

\begin{proof}
W.l.o.g., we only consider quasi-normal derivations. Initial sequents other than regularity and truth rules are derivable in a finite number of steps via the fixed-point rules. For instance, suppose the derivation ends with 
        \begin{IEEEeqnarray*}{L}
	\def\fCenter{ \mbox{$\Gamma$} }
		\ax{ \fCenter, \Tz\val s}
			\un{ \fCenter, \Tz\uTz s }
		\Dis
        \end{IEEEeqnarray*}
By IH we have $\h{2}\svdash{\star}{<\omega\cdot k}\Gamman, \p0{}\val s$. By form of $\Acr_0$, this yields $\h{2}\svdash{\star}{<\omega\cdot k+c}\Gamman, \Acr_0(\uTz s, \p0{})$, and hence $\h{2}\svdash{\star}{<\omega\cdot k+c+1}\Gamman, \p0{}\uTz s$ by the fixed-point rule. The cut rule is immediate by IH, and $\forall$ becomes an application of the $\omega$-rule: since the height $<\omega\cdot k$ is uniform for all $n$, this adds only one to the height. 

The height for the derivability of regularity depends on the complexity of \vphi but can be bounded by $\omega + c$. As for internal induction, suppose the derivation ends with
    \begin{IEEEeqnarray*}{C}
        \axc{\svdash{}{m}\Gamma, \Fml^1(\vphi)}
        \axc{\svdash{}{m}\Gamma, \Tu\vphi(0)}
        \axc{\svdash{}{m}\Gamma,\neg\Tu\vphi(y), \Tu\vphi(y+1)}
        \trc{\svdash{}{h}\Gamma, \Tu\vphi(r)}
        \Dis
    \end{IEEEeqnarray*}
By IH we have, for a suitably large $k$, $\h{2}\svdash{\star}{<\omega\cdot(k+1)}\Gamman, \p1{}\vphi(\ovl 0)$ and $\h{2}\svdash{\star}{<\omega\cdot(k+1)}\Gamman, \neg\p1{}\vphi(\ovl n), \p1{}\vphi(\ovl{n+1})$ for all $n$. Now let $r$ be an arbitrary term with value $p$. The desired conclusion can be obtained by $p$-many cuts on $\p1{}$ atoms against the premises.
\end{proof}

\begin{crl}\label{crl:embedding_kftint-h2}
    If $\kftintcs\vdash\sigma$ for $\sigma\in\lpa$, then $\h2\svdash{\star}{<\omega\cdot k}\neg\Const,\sigma$ for some $k<\omega$.
\end{crl}

Next, we asymmetrically embed \h2 into a semi-formal system \ah2 (cf. \cite{can85Note, fuj11Autonomous}). The Tait-style language $\Lcr_{\ah2}$ of \ah{2} is obtained by closing off $\lpa \cup \{\p0{}\}\cup\{\p1\eta\mid\eta<\omega^\omega\}$ under logical operations.
Axioms and rules of \ah2 are the following:

\noindent
\textit{Axioms.}
	\begin{enumerate}[label=Ax\arabic*]
	\setcounter{enumi}{-1}
		\item $\Gamma, \vphi$ for a true \lpa-literal $\vphi$.
		\item $\Gamma, \neg\p0{}t, \p0{}s$ for $s^\Nbl=t^\Nbl$.
        \item $\Gamma, \neg\p10t$ for all closed terms $t$.
	\end{enumerate}
\noindent
\textit{Logical rules and Cut:} Same as those of $\omega$-arithmetic.

\noindent
\textit{Fixed-point rules for \p0{}}: Same as those of \h2.

\noindent
\textit{Fixed-point rules for $\p1{\eta}$}:
	\begin{IEEEeqnarray*}{L}
	\def\fCenter{ \mbox{$\Gamma$} }
		\ax{ \fCenter, \Acr_1(s, \p0{}, \p1{\eta})}
		\rl{\msf{FP_1^+}}
			\un{ \fCenter, \p1{\eta+1}s }
		\Dis
	\hspace*{20mm}
	\def\fCenter{ \mbox{$\Gamma$} }
		\ax{ \fCenter, \neg\Acr_1(s, \p0{}, \p1{\eta})}
        \rl{\msf{FP_1^-}}
			\un{ \fCenter, \neg\p1{\eta+1}s }
		\Dis
	\end{IEEEeqnarray*}

\noindent
\textit{Limit rules for $\p{1}{\lambda}$, for $\lambda$ limit}:
	\begin{IEEEeqnarray*}{L}
	\def\fCenter{ \mbox{$\Gamma$} }
		\ax{ \fCenter, \p{1}{\eta}s \hspace*{3mm} \text{for some } \eta<\lambda}
        \rl{\msf{Lim_1^+}}
			\un{ \fCenter, \p{1}{\lambda}s }
		\Dis
	\hspace*{20mm}
	\def\fCenter{ \mbox{$\Gamma$} }
		\ax{ \fCenter, \neg\p{1}{\eta}s \hspace*{3mm} \text{for all } \eta<\lambda}
        \rl{\msf{Lim_1^-}}
			\un{ \fCenter, \neg\p{1}{\lambda}s }
		\Dis
	\end{IEEEeqnarray*}
    
For each $\Lcr_{\ah2}$-formula \vphi, we define its $\msf{A}$-rank, $\mis{rk}_\Asf(\vphi)$, as follows: if \vphi is a literal of the form $\p{1}\eta(t)$ or $\neg \p{1}\eta(t)$, then $\mis{rk}_\Asf(\vphi)=\omega\cdot\eta$; if \vphi is a literal of any other form, $\mis{rk}_\Asf(\vphi)=0$; for a compositional formula \vphi, $\mis{rk}_\Asf(\vphi)$ is inductively defined in the standard way. We write $S\svdash{\rho}{\delta}\Gamma$ (respectively $S\svdash{\rho}{<\delta}\Gamma$) when $\Gamma$ is derivable in $S$ with a derivation of height $\leq\delta$ (resp.~$<\delta$) and such that every cut formula has \Asf-rank $<\rho$. Predicative cut-elimination (see e.g. \cite[§22]{schu77Proof}) removes all cuts of \Asf-rank $\geq1$, since cuts on $\p1{\eta+1}s$ or $\p1{\lambda}s$ reduce to cuts of lower \Asf-rank via \msf{FP_1^\pm} and \msf{Lim_1^\pm}, respectively. Cuts on $\p10$-literals and false arithmetical literals can be deleted as they are never principal. Since \p0{} is not stratified, cuts on $\p0{}s$ cannot be reduced. For derivations in which cuts are restricted to literals of the form \p0{}, we write $\ah2\svdash{\star}{\delta}\Gamma$ and we say that the derivation is quasi-normal.

Given a formula $\vphi\in\Lcr_{\h2}$, let $\vphi[\beta, \gamma]$ be the $\Lcr_{\ah2}$-formula obtained from \vphi by replacing occurrences of \p1{} as follows:
	\begin{IEEEeqnarray*}{RCL+RCL}
		\neg\p1{}(s) &\mapsto& \neg\p1\beta(s),
		& \p1{}(s) &\mapsto & \p1\gamma(s).
	\end{IEEEeqnarray*}

In order to asymmetrically embed \h2 into \ah2 via this mapping, we need to prove the so-called \textit{Persistence Lemma}. We begin with an auxiliary lemma.

\begin{lemma}\label{le_aux-persistence}
The following hold for all $s$ and all $\vphi\in\Lcr_{\h2}$:
\begin{enumeratei}
    \item If $\zeta\leq\eta$, then $\ah2\svdash{0}{<\omega\cdot(\zeta+1)}\neg\p1{\zeta}(s), \p1{\eta}(t)$ for all closed terms $s,t$ with $s^\Nbl= t^\Nbl$;
	\item Assume $\beta\leq\beta'\leq\gamma'\leq\gamma$. Then $\ah2\svdash{0}{<\omega\cdot(\gamma+1)}\neg\vphi[\beta', \gamma'], \vphi[\beta, \gamma]$.
\end{enumeratei}
\end{lemma}

\begin{proof}
    (i) By induction on $\zeta$. For $\zeta=0$, this is just an instance of Ax2. If $\zeta=\zeta'+1,\eta=\eta'+1$, by IH we have the claim for $\zeta'$. By positivity of \p1{} in $\Acr_1$ and monotonicity, we thus obtain a cut-free derivation of $\neg\Acr_1(s, \p0{},\p1{\zeta'}), \Acr_1(t, \p0{},\p1{\eta'})$ for all $\eta'\geq\zeta'$. The desired conclusion follows by \msf{FP_1^-} and \msf{FP_1^+}. If $\zeta=\zeta'+1$ and $\eta$ limit, then $\zeta<\eta$ and we apply \msf{Lim^+_1} to $\eta:=\zeta$. If $\zeta$ is a limit, apply \msf{Lim^-_1} to all $\xi<\zeta$, whose heights are $<\omega\cdot(\xi+1)\leq\omega\cdot\zeta$. 
    
    (ii) Straightforward induction on $\vphi$, using (i).
\end{proof}

For $\Gamma=\varphi_1,\dots,\varphi_m$ and $\vec\beta=(\beta_1,\dots,\beta_m)$, $\vec\gamma=(\gamma_1,\dots,\gamma_m)$ sequences of ordinals, let
    \[
        \Gamma[\vec\beta,\vec\gamma]:=\varphi_1[\beta_1,\gamma_1],\dots,\varphi_m[\beta_m,\gamma_m]. 
    \]
Inequalities involving sequences are meant componentwise, and a single ordinal stands for the constant sequence; so $\Gamma[\zeta,\eta]$ is the special case of constant sequences.

\begin{lemma}[Persistence]\label{le_persistence}
Let $\vec\beta\le\vec\beta'\le\vec\gamma'\le\vec\gamma\le\eta$ and suppose $\ah2\svdash\rho{\delta}\Gamma[\vec\beta',\vec\gamma'],\Delta$. Then
$\ah2\svdash{\max\{\rho,\,\omega\cdot(\eta+1)\}}{<\max\{\delta+\omega,\,\omega\cdot(\eta+1)\}}\Gamma[\vec\beta,\vec\gamma],\Delta$.
\end{lemma}

\begin{proof}
For $i=1,\dots,m$, cut $\varphi_i[\beta_i',\gamma_i']$ against the sequent $\neg\varphi_i[\beta_i',\gamma_i'],\varphi_i[\beta_i,\gamma_i]$ of \Cref{le_aux-persistence}(ii), which has height $<\omega\cdot(\eta+1)$. The cut formulae have \Asf-rank at most $\omega\cdot\gamma_i'+c(\varphi_i)<\omega\cdot(\eta+1)$, where $c(\vphi)<\omega$ is the surface complexity of $\vphi$,
and the $m$ cuts add only finitely many steps (cf.~\cite[Lemma 3.7]{can85Note} for a similar argument).
\end{proof}

\begin{theorem}[Asymmetric interpretation]\label{th_asym}
    If $\h{2}\svdash\star{\delta}\Gamma$ with $\delta<\omega^2$, then, for each $\zeta<\omega^\omega$, $\ah{2}\svdash{\omega\cdot(\zeta+2^\delta+1+\delta)}{\omega\cdot(\zeta+2^\delta+1+\delta)}\Gamma[\zeta, \zeta+2^\delta]$.
\end{theorem}

\begin{proof}[Proof Sketch]
By induction on $\delta$. Write $\theta(\zeta,\delta):=\omega\cdot(\zeta+2^\delta+1+\delta)$. Since $\delta=\omega\cdot j+n$ for some $j,n<\omega$, we have $2^\delta=\omega^j\cdot2^n$, so $\zeta+2^\delta<\omega^\omega$ (as $\omega^\omega$ is closed under addition). Hence $\Gamma[\zeta,\zeta+2^\delta]$ is an $\Lcr_{\ah2}$-sequent, and the IH is only applied to $\delta'<\delta$ and $\zeta'\le\zeta+2^\delta$.
For the base step where $\Gamma$ is an axiom, it suffices to consider the case where $\Gamma=\Gamma',\neg\p1{}s, \p1{}t$ for $s^\Nbl=t^\Nbl$. The desired conclusion $\Gamma'[\zeta, \zeta+1],\neg\p1{\zeta}s, \p1{\zeta+1}t$ then follows from \Cref{le_aux-persistence}.

For $\delta>0$, we rely on the following consequence of the Persistence Lemma:
    \begin{IEEEeqnarray}{L}
        \text{Let } \delta'<\delta, \zeta\le\zeta' \text{ and } \zeta'+2^{\delta'}\le\zeta+2^{\delta}; 
        \text{ and let  }\zeta\le\vec\beta\le\vec\gamma\le\zeta+2^{\delta}. \nonumber
        \\
        \text{If } \ah2\svdash{\theta(\zeta',\delta')}{<\theta(\zeta',\delta')+\omega}\Gamma[\vec\beta,\vec\gamma],\Delta, 
        \text{ then }\ah2\svdash{\theta(\zeta,\delta)}{<\theta(\zeta,\delta)}\Gamma[\zeta,\zeta+2^{\delta}],\Delta. \label{eq_asymcut}
    \end{IEEEeqnarray}
Its proof is \Cref{le_persistence} with $\eta:=\zeta+2^\delta$, using $\theta(\zeta',\delta')+\omega,\omega\cdot(\zeta+2^\delta+1)\leq\theta(\zeta,\delta)$.

The argument is then standard. We only show two cases as an example. First, suppose that the derivation ends with the positive fixed-point rule for \p1{}:
    \begin{IEEEeqnarray*}{L}
        \def\fCenter{ \mbox{ $\h2$ } }
		\ax{ \fCenter\svdash{\star}{\delta'}\Gamma, \Acr_1(s, \p0{},\p1{}) }
			\un{ \fCenter\svdash{\star}{\delta}\Gamma, \p1{}(s) }
		\Dis
    \end{IEEEeqnarray*}
Since \p1{} occurs only positively in $\Acr_1(x,\p0{},\p1{})$, by IH we have
    \[
        \ah2\svdash{\theta(\zeta,\delta')}{\theta(\zeta,\delta')}\Gamma[\zeta, \zeta+2^{\delta'}], \Acr_1(s, \p0{}, \p1{\zeta+2^{\delta'}}).
    \]
By $\msf{FP_1^+}$
    \[
        \ah2\svdash{\theta(\zeta,\delta')}{\theta(\zeta,\delta')+1}\Gamma[\zeta, \zeta+2^{\delta'}], \p1{\zeta+2^{\delta'}+1}(s).
    \]
Since $2^{\delta'}+1\leq2^\delta$, the conclusion follows by \eqref{eq_asymcut}.

Second, suppose the derivation ends with Cut on \p1{}:
    \begin{IEEEeqnarray*}{L}
        \def\fCenter{ \mbox{ $\h2$ } }
		\ax{ \fCenter\svdash{\star}{\delta_0}\Gamma, \p1{}(s) }
		\ax{ \fCenter\svdash{\star}{\delta_1}\Gamma, \neg\p1{}(s) }
			\bn{ \fCenter\svdash{\star}{\delta}\Gamma }
		\Dis
    \end{IEEEeqnarray*}
for $\delta_0,\delta_1<\delta$. By IH, we have
    \begin{IEEEeqnarray*}{L}
        \ah2\svdash{\theta(\zeta,\delta_0)}{\theta(\zeta,\delta_0)}\Gamma[\zeta, \zeta+2^{\delta_0}], \p1{\zeta+2^{\delta_0}}(s),
        \\
        \ah2\svdash{\theta(\zeta+2^{\delta_0},\delta_1)}{\theta(\zeta+2^{\delta_0},\delta_1)}\Gamma[\zeta+2^{\delta_0}, \zeta+2^{\delta_0}+2^{\delta_1}], \neg\p1{\zeta+2^{\delta_0}}(s).
    \end{IEEEeqnarray*}
By \eqref{eq_asymcut}, we obtain
    \begin{IEEEeqnarray*}{L}
        \ah2\svdash{\theta(\zeta,\delta)}{<\theta(\zeta,\delta)}\Gamma[\zeta,\zeta+2^{\delta}],\p1{\zeta+2^{\delta_0}}(s)
        \\
        \ah2\svdash{\theta(\zeta,\delta)}{<\theta(\zeta,\delta)}\Gamma[\zeta,\zeta+2^{\delta}],\neg\p1{\zeta+2^{\delta_0}}(s).
    \end{IEEEeqnarray*}
The desired conclusion follows by a cut of \Asf-rank $\omega\cdot(\zeta+2^{\delta_0})<\theta(\zeta,\delta)$.
\end{proof}

Having shown how to asymmetrically embed \h2 into \ah2 with the given bounds, we now adapt the conservativity argument by \cite{fuj11Autonomous}.

\begin{lemma}[{\cite[Lemma 8]{fuj11Autonomous}}]\label{le:Fuj8}
    For $\rho>0$, for $\eta_0,\eta_1$ arbitrary, if (a) $\ah2\svdash{\rho}{\delta_0}\Gamma, \p1{\eta_0}s$ and (b) $\ah2\svdash{\rho}{\delta_1}\Gamma, \p1{\eta_1}\ng s$, then $\ah2\svdash{\rho}{\delta_0\#\delta_1}\Gamma$ (for $\delta_0\#\delta_1$ the natural sum of $\delta_0$ and $\delta_1$).
\end{lemma}

\begin{lemma}\label{le:Cons2}
    For $\eta<\omega^m$ with $1\leq m<\omega$: $\ah2\svdash{1}{<\omega^{m+1}}\forall x\in\Sent_{\lkft}\left(\neg\p1{\eta}x\lor\neg\p1{\eta}\ng x\right)$.
\end{lemma}

\begin{proof}
    By \Cref{le_aux-persistence}(i), we have, for all $s$, $\ah2\svdash{1}{\delta_0}\neg\p1\eta s, \p1\eta s$ and $\svdash{1}{\delta_1}\neg\p1\eta\ng s, \p1\eta\ng s$ for $\delta_0, \delta_1<\omega\cdot(\eta+1)$. By \Cref{le:Fuj8}
    these yield $\svdash{1}{\delta_0\#\delta_1}\neg\p1\eta s, \neg\p1\eta\ng s$, with $\delta_0\#\delta_1<\omega^{m+1}$. The desired conclusion follows by the rules $\lor$ and $\omega$.
\end{proof}

Let $\Consu:=\forall x\in\Sent_{\lkf}(\neg\Tz x\lor\neg\Tz \ng x)$. Then, writing \Consu instead of $\Consu[\p0{}/\Tz]$:

\begin{crl}[Cf.~{\cite[Lemma 10]{fuj11Autonomous}}]\label{crl:Fuj10}
    Let $1\leq k<\omega$ and $\delta<\omega\cdot k$, and let $\Gamma$ be a set of sentences not containing \p1{}. If $\h2\svdash{\star}{\delta}\Gamma,\neg\Const$, then $\ah2\svdash{\omega^{k+1}}{<\omega^{k+1}}\Gamma,\neg\Consu$.
\end{crl}

\begin{proof}
Let $\delta = \omega\cdot j+n$ with $j<k$ and recall that $\theta(\zeta,\delta):=\omega\cdot(\zeta+2^\delta+1+\delta)$. Then $2^\delta=\omega^j\cdot2^n<\omega^k$ and $\delta<\omega^k$. Since $2^\delta+1+\delta<\omega^k$, we have $\theta(0,\delta)<\omega^{k+1}$. Choosing $\zeta:=0$ in \Cref{th_asym}, from the assumption we get 
    \[
        \ah2\svdash{\theta(0,\delta)}{\theta(0,\delta)}\Gamma,\neg\Const[0,2^\delta].
    \]
By $\lor$-inversion, 
    \[
        \ah2\svdash{\theta(0,\delta)}{\theta(0,\delta)}\Gamma,\neg\Consu, \exists x\in\Sent_{\lkft}\left( \p1{2^\delta}x\land \p1{2^\delta}\ng x\right).
    \]
By \Cref{le:Cons2} with $m=k$, 
    \[
        \ah2\svdash{1}{<\omega^{k+1}}\forall x\in\Sent_{\lkft}\left( \neg\p1{2^\delta}x\lor \neg\p1{2^\delta}\ng x\right).
    \]
A cut of \Asf-rank $\omega\cdot2^\delta+3<\omega^{k+1}$ yields the desired conclusion.
\end{proof}

\begin{lemma}[Cf.~{\cite[Lemma 11]{fuj11Autonomous}}]\label{le:h2-to-h1}
    Let $1\leq k<\omega$ and $\delta<\omega\cdot k$, and let $\Gamma$ be a set of sentences not containing \p1{}. If $\h2\svdash{\star}{\delta}\Gamma,\neg\Const$, then $\h1\svdash{\star}{<\vphi_{k+2}0}\Gamma, \neg\Consu$.
\end{lemma}

\begin{proof}
    If $\h2\svdash{\star}{\delta}\Gamma,\neg\Const$, then $\ah2\svdash{\omega^{k+1}}{<\omega^{k+1}}\Gamma,\neg\Consu$ by \Cref{crl:Fuj10}. By Tait's refinement of cut-elimination (see \cite[p.~204]{schu77Proof}), $\ah2\svdash{1}{<\vphi_{k+1}(\omega^{k+1})}\Gamma,\neg\Consu$. Since any cut of \Asf-rank 0 except those on \p0{}-atoms can be removed, we get $\ah2\svdash{\star}{<\vphi_{k+1}(\omega^{k+1})}\Gamma,\neg\Consu$.
    Since $\Gamma,\neg\Consu\subseteq\Lmc_{\h1}$, by quasi-normality we can replace \ah2 with \h1, so we are done by $\vphi_{k+1}(\omega^{k+1}) < \vphi_{k+2}0$.
\end{proof}

We can finally prove the upper bound theorem:

\begin{proof}[Proof of \Cref{th:UB-kftintcs}]

    If $\kftintcs\vdash\sigma$, then $\h2\svdash{\star}{<\omega\cdot k}\neg\Const, \sigma$ by \Cref{crl:embedding_kftint-h2}. By \Cref{le:h2-to-h1}, we infer $\h1\svdash{\star}{\gamma}\neg\Consu, \sigma$ for $\gamma<\vphi_{k+2}0$. By \cite[Lemmata 11, 15]{fuj11Autonomous}, for $\varepsilon(\gamma)$ the least $\varepsilon$-number above $\gamma$ we obtain $\h0\svdash{\star}{<\vphi_{\varepsilon(\gamma)}0}\sigma$. Since $\varepsilon(\gamma)\leq\vphi_{k+2}0<\vphi_\omega0$, we have $\vphi_{\varepsilon(\gamma)}0<\vphi_{\vphi_\omega0}0$. Since all transformations are primitive recursive on codes of derivations (cf.~\cite[p.~926]{fuj11Autonomous}), it follows that $\PA+\TI(<\vphi_{\vphi_\omega0}0)\vdash\sigma$.
\end{proof}

This yields:

\begin{theorem}
    $\vert\kftint\vert = \vert\kftintcs\vert = \vert\kftintcssp\vert = \vert\kfdint\vert = \vert\kfdintcs\vert = \vphi_{\vphi_\omega0}0$.
\end{theorem}

We expect that the above interpretability results can be used to improve the lower bound of Field's \INT. In fact, while Field (cf. \cite[p.~257]{fie22Power}) conjectures that \INT is arithmetically equivalent to \kf, whose ordinal is $\vphi_{\varepsilon_0}0$, we conjecture that $\vert\INT\vert=\vert\kftint\vert$.\footnote{
    Carlo Nicolai has informed me that he has independently arrived at a similar conjecture. We will deepen the analysis in forthcoming joint work.
}

\section{Concluding remarks and further developments}

The notion of determinateness has typically been assumed to be \textit{definable} within \kfcs as truth or falsity. And it is not hard to see that each of the determinateness axioms from \Cref{dfn:kfd} but the right conjunct of \msf{D8} is derivable in \kfcs under the definition $\D x:\lra\T x\lor\F x$. Yet, crucially, no sentence \vphi can be proven to be \textit{truly non-determinate} in \kfcs. That is, there is no \vphi such that $\kfcs\vdash\T\neg\D\vphi$. And there could be no such \vphi, since $\kfcs\vdash\T(\neg\T\vphi\land\neg\F\vphi)\ra\F\vphi\land\T\vphi$.

This means that the equivalence between $\F\D\vphi$ and $\neg\D\vphi$ is rather crucial, from both a conceptual and a proof-theoretical point of view. Conceptually, a key aspect of \kfd is that there are two ways of being true: sentences can be \textit{true}, or they can be \textit{true-and-determinate}. Importantly, the principle $\F\D\vphi\lra\neg\D\vphi$ yields that the predicate `true' properly extends the predicate `true-and-determinate'. Specifically, for any \vphi such that $\neg\D\vphi$, we have $\T\neg\D\vphi\land\neg\D\neg\D\vphi$. The right conjunct is a consequence of \msf{D3} and \msf{D4}, which from $\D\neg\D\vphi$ yield $\D\vphi$.
This conceptual difference is clearly displayed in the proof-theoretic analysis. In fact, the principle $\F\D\vphi\lra\neg\D\vphi$ guarantees that \textit{true} is fully compositional and disquotational over \textit{true-and-determinate}. Technically, this translates into the fact that \kfd can interpret \kft, by mapping \Tz to \textit{true-and-determinate} and \Tu to \textit{true}. 
The same considerations apply also to \cdp and \cd: only the former derives the equivalence between $\F\D\vphi$ and $\neg\D\vphi$.
Accordingly, while \cd interprets \kf, \cdp interprets $\CT[\kfcs]$ \cite{fh24Classical}. Of course, also in \cd one can distinguish between \T and $\T\cap\D$. However, the former is not a Tarskian truth predicate for $\T\cap\D$: given the lack of the equivalence between $\F\D\vphi$ and $\neg\D\vphi$, there are models of \cd in which, e.g., for a liar sentence $\lambda$, $\T\D\lambda$ holds but $\D\lambda$ fails.

\subsection{Further directions}
Let us conclude with some possible extensions of the present study.
    \begin{enumeratei}
        \item \emph{Minimality.} The characterization of models of \kfd from §\ref{sec:semantics} suggests that the theory can be expanded with two different minimality schemata, namely one for \D and one for \T. What is the strength of the resulting theories?
        \item \emph{Bi-interpretability.}
        There is a clear sense in which adding \D to a truth theory amounts to adding a new truth predicate. Just how strong is this sense in the case of \kfdcs? Specifically, are \kfdcs---possibly extended with some persistence principle---and \kftcssp bi-interpretable? 
        \item \emph{Alternative schemes.} What are natural axiomatizations of determinateness over \wkf \cite{fuj10Relative} and \msf{VF} \cite{can90Theory}, respectively, and, what is the strength of the resulting theories? Note that, while the limitative results by \cite{cgw26} entail that the theory obtained from \kfd by replacing \msf{D3} with $\forall s\D\ud\D s$ is inconsistent, they do not entail that the same holds for \wkf. Thus, it is of interest to consider an axiomatization of  determinateness over \wkf with the axiom $\forall s\D\ud\D s$.
    \end{enumeratei}

\bigskip
\subsection*{Acknowledgments}
    [TO BE ADDED]

\bigskip
\bibliography{bib}
\bibliographystyle{alpha}
\end{document}